\pdfoutput=1
\documentclass[a4paper,11pt]{article}
\usepackage{amsmath,amssymb,amsthm,graphicx}
\usepackage{hyperref,xurl,enumitem,bm,esint}
\usepackage{color}
\usepackage[dvipsnames]{xcolor} %
\usepackage{comment}
\usepackage{algorithm}
\usepackage{algpseudocode}
\usepackage[T1]{fontenc}
\usepackage[normalem]{ulem}
\usepackage{pifont}
\usepackage{xcolor,enumitem}
\usepackage{graphicx}
\usepackage{multirow}
\usepackage{makecell}
\usepackage{subcaption}
\usepackage{mwe}
\usepackage{overpic}
\usepackage[margin=30mm]{geometry}
\usepackage{booktabs}
\usepackage{placeins}
\usepackage{pgfplots}               %
\usepgfplotslibrary{groupplots}     %
\pgfplotsset{compat=1.16}
\definecolor{stablecolor}{RGB}{204,232,243}
\definecolor{oscillatorycolor}{RGB}{252,216,216}

\newcommand{\tnorm}[1]{|\!|\!| #1 |\!|\!|}

\newcommand{\bd}{\mathbf d}

\newcommand{\bb}{\mathbf b}

\newcommand{\bx}{\bm{x}}
\newcommand{\bn}{\mathbf n}
\newcommand{\br}{\mathbf r}
\newcommand{\bA}{\mathbf A}
\newcommand{\bp}{\mathbf p}
\newcommand{\bH}{\mathbf H}
\newcommand{\Ast}{A^{\mathrm{st}}}

\usepackage{todonotes}
\theoremstyle{plain}

\newtheorem{lem}{Lemma}[section]
\newtheorem{prop}{Proposition}[section]

\newtheorem{remark}{Remark}[section]
\newtheorem{assump}{Assumption}[section]

\numberwithin{equation}{section}

\def\cl {\nonumber \\}

\definecolor{diff}{HTML}{C62828}

\title{Coupled and fully decoupled
approaches for a chemotaxis-type model of urban crime}

\author{Tim B\"{u}rchner$^1$, Baoli Hao$^2$, Yunhui He$^3$,  Ming Zhong$^3$, and Annalisa Quaini$^1$}

\date{}
\begin{document}
\maketitle
\pagestyle{myheadings}
\markright{Coupled and decoupled iterations for an urban crime model}

\begin{center}
\noindent $^1$ Chair of Computational Methods in Civil and Environmental Engineering, Technical University of Munich,
Arcisstra{\ss}e 21, 80333 Munich, Germany \\ {\tt tim.buerchner@tum.de; annalisa.quaini@tum.de}

\vskip .2cm
$^2$ Department of Applied Mathematics, Illinois Institute of Technology, 10 West 35th St,
Chicago, IL 60616, USA \\ {\tt bhao2@hawk.illinoistech.edu}

\vskip .2cm
\noindent $^3$ Department of Mathematics, University of Houston, 3551 Cullen Blvd, Houston, TX 77204, USA \\
{\tt  yhe43@central.uh.edu; mzhong3@central.uh.edu}
\end{center}

\begin{abstract}
We study the iterative solution of the nonlinear systems arising at each time step from the finite element discretization of a
model of urban crime with police.
The model couples three PDEs for attractiveness, criminal density, and police density to an ODE for the reported crime rate.
We compare a fully decoupled fixed-point iteration, which solves for one variable at a time, with a coupled iteration, in which attractiveness and criminal density are computed together through a Newton linearization of the reaction term.
Under suitable bounds on the discrete solution and for sufficiently small $\Delta t$, we prove that both iterations converge with an asymptotic rate of order $\Delta t$.
To improve robustness, we combine both iterations with Anderson acceleration (AA), with the least-squares problem posed in a block-weighted norm.
Numerical results for ten parameter sets  confirm the predicted dependence of the convergence rate on $\Delta t$.
Without AA, the fully decoupled iteration fails in one case, in which the coupled iteration converges. With AA, both approaches complete all the simulations. Our numerical tests show that
the fully decoupled approach with AA is the most efficient.
\end{abstract}


\vskip .2cm
\noindent \emph{AMSC}: 65H10, 65M60, 35Q91

\vskip .5cm
\centerline{Dedicated to the memory of Yuri A.~Kuznetsov}
\vskip .5cm

\section{Introduction}\label{sec:intro}
This paper is dedicated to the memory of Yuri A.~Kuznetsov. His monograph with Gury Marchuk on iterative methods and quadratic functionals~\cite{MarchukKuznetsov1972} is an early expression of a theme that runs through his work: the choice of an iteration should reflect the structure of the discrete problem. His work on the approximate solution of implicit difference schemes~\cite{Kuznetsov1988}, on algebraic multigrid and domain decomposition~\cite{Kuznetsov1989amg,Kuznetsov1989cmame}, and on fictitious domain methods~\cite{HeikkolaEtAl1998} gave concrete form to this principle for large sparse systems. A question in that tradition motivates the present study. When several fields are coupled in an implicit time step, which of them should be solved for together, and which couplings can be left to an outer iteration? We examine this question for a model of urban crime in which the choice of splitting changes both the linear systems to be solved and the convergence of the nonlinear iteration.

In this paper, we consider the PDE model of urban crime with police introduced in~\cite{Hao2026b}, which extends the continuum model of residential burglary in~\cite{short2008statistical}.
The latter couples the density of criminals to the attractiveness of residential sites: criminals move preferentially toward more attractive sites, and burglaries increase local attractiveness.
This feedback is of the same type as in chemotaxis models~\cite{KELLER1970399}, where a population moves up the gradient of a signal that it helps to generate.
The model in~\cite{short2008statistical} can lead to crime hotspots, i.e., small regions where criminals concentrate and crime persists over time.
This model has been studied in terms of hotspot formation and suppression~\cite{short2010nonlinear,short2010dissipation} and local well-posedness~\cite{rodriguez2010}, and extended to complex geometries~\cite{Hao2026}.
Law enforcement has been incorporated into both agent-based and continuum models~\cite{jones2010,PITCHER_2010,zipkin2014,Hao2026b}.

The model in~\cite{Hao2026b} adds police as a second population that responds to a different signal.
While criminals follow the gradient of the logarithm of attractiveness, police follow the gradient of the logarithm of the reported crime rate, and their presence deters burglary.
The reported rate relaxes toward the instantaneous crime rate on a prescribed time scale, so the signal guiding police can persist after local crime has subsided.
The result is a system of three partial differential equations for attractiveness, criminal density, and police density coupled to an ordinary differential equation for the reported crime rate.
In an implicit time step, the change in each population is coupled to the change in the field that directs its motion.
This coupling is most delicate near a sharp hotspot, where the logarithmic gradients that guide the two populations are large and highly sensitive to small changes in the signals.

The numerical literature on Keller--Segel systems offers several ways to handle this coupling.
Conservative upwind and flux-corrected finite element methods address transport and the preservation of nonnegative densities~\cite{Saito2007,Strehl2010}.
Of particular relevance here is~\cite{StrehlSokolovTurek2012}, which compares monolithic Newton-like and decoupled finite element approaches for chemotaxis
and examines the cost and stability of solving the population and signal equations together.
Other works exploit the structure of the algebraic systems: in~\cite{HuZhang2023}, the density update of a parabolic-elliptic Keller--Segel discretization is reformulated as a symmetric positive definite system
and solved by preconditioned conjugate gradients, while~\cite{ZhangGeYang2023} develops a multigrid algorithm for a fourth-order compact difference discretization of the time-dependent Keller--Segel system.
In the related setting of chemotaxis optimal control, block preconditioners for GMRES applied to the saddle-point systems arising in a Gauss--Newton iteration are constructed in~\cite{DolgovPearson2019}.
We focus instead on the convergence of the nonlinear iteration for the full four-field model.

Our starting point is the partitioned scheme used in~\cite{Hao2026b}, which extends the finite element approach of~\cite{Hao2026}. After backward Euler discretization in time, each iteration updates attractiveness, criminal density, reported crime, and police density in succession. The reported crime rate is updated at the nodes, so each iteration requires three scalar linear solves. Only the criminal and police densities need to be passed from one iteration to the next. We compare this fully decoupled iteration with a second approach in which a Newton linearization of the reaction product couples attractiveness and criminal density in a single block solve. We refer to it as the coupled approach. At each time step, both approaches converge to the same discrete solution. The two approaches differ in the attractiveness that drives the motion of criminals toward attractive sites: the fully decoupled approach uses the attractiveness just computed in the same iteration, whereas the coupled approach uses the attractiveness from the previous iteration.

At difficult time steps, e.g., in the presence of sharp hotspots, both iterative approaches may converge slowly or even fail.
To improve their robustness at a low cost, we combine both approaches with Anderson acceleration (AA)~\cite{Anderson1965,Walker2011}.
AA treats the fixed-point map as a black box and computes each new iterate from a short history of previous iterates by solving a small least-squares problem,
so it requires no changes to the linear solvers. It has been used to accelerate nonlinear fixed-point iterations and splitting schemes in several applications,
e.g., the Navier--Stokes equations~\cite{Pollock2019}, variably saturated flow~\cite{Lott2012}, and poromechanics~\cite{Both2019}.
The local convergence theory of AA assumes that the underlying map is contractive~\cite{TothKelley2015}
and relates the improvement in the convergence rate to the reduction achieved in the least-squares problem~\cite{Evans2020}.
This motivates the convergence analysis of the two approaches presented in this paper.
Since the unknowns of our problem have very different magnitudes, we pose the least-squares problem in a block-weighted norm that scales each field by its size at the previous time step.

Under suitable bounds on the discrete solution and for sufficiently small $\Delta t$, we prove that both decoupled and coupled iterations converge with a rate proportional to the time step $\Delta t$.
For the fully decoupled approach, the error is reduced by a factor of order $\Delta t$ at every iteration.
For the coupled approach, this holds only asymptotically: since the equation for the criminal density uses the attractiveness from the previous iteration,
the first iteration may reduce the error less, but the asymptotic rate is again of order $\Delta t$.

We test both approaches on ten parameter sets that produce stationary and oscillatory dynamics.
Without acceleration, most time steps require only a few iterations. However,
in one oscillatory case, the fully decoupled iteration fails to converge at a time step where the coupled iteration converges.
With AA, the fully decoupled approach converges at that time step as well.
Hence, the main benefit of AA in our tests is robustness: it restores convergence at difficult time steps, while it reduces the average number of iterations only moderately.
Overall, the fully decoupled approach with AA is the most efficient in our tests.

The rest of the paper is organized as follows.
Sec.~\ref{sec:pd} presents the model and its discretization. Sec.~\ref{sec:AA} introduces the acceleration procedure. Secs.~\ref{sec:decoupled} and~\ref{sec:coupled} derive the two fixed-point maps, describe their linear systems, and analyze convergence. Sec.~\ref{sec:num_res} presents numerical results that support the analysis in Secs.~\ref{sec:decoupled} and~\ref{sec:coupled}. Sec.~\ref{sec:concl} summarizes the findings.
\section{Problem definition}\label{sec:pd}

This section briefly recalls the PDE model from \cite{Hao2026b}.
We consider a bounded, connected Lipschitz domain $\Omega\subset \mathbb{R}^2$ and a time interval  $(0,T)$.
Let $A(t, \bx)>0$ denote the level of attractiveness. This variable, which measures the criminals’
perception at time $t$ of the attractiveness of a site located at $\bx$, is the sum of a static component $\Ast >0$ and a dynamic component associated with repeat
and near-repeat victimization. Let
$\rho(t, \bx)$ and $\pi(t, \bx)$ be the density
of criminals and police officers, respectively.
Let $H(t, \bx)>0$ be the number of crimes per unit area per unit time lagged in time.
The PDE model in \cite{Hao2026b} was derived from the mean-field limit of an agent-based model, which assumes that if police are not busy at a certain location, they move according to
a biased random walk towards areas with higher reported crime rates. Variable $H$ is meant to account for the fact that crime data become available and are used for decision-making with some delay, denoted by $\tau$.
All variables are  nondimensionalized.
The nondimensional PDE model from \cite{Hao2026b} reads:
find $A(t, \bx)$, $\rho(t, \bx)$, $\pi(t, \bx)$, and $H(t, \bx)$, such that
\begin{align}
     & \frac{\partial A}{\partial t} - \eta \Delta A + A - \rho A e^{-\pi}= - \eta \Delta \Ast + \Ast, \label{eq:continous_1}\\
     &\frac{\partial \rho}{\partial t} - \nabla \cdot \Big(\nabla\rho - \frac{2 \nabla A}{A} \rho\Big)+\rho A e^{-\pi}  - q e^{-\pi}= 0,\label{eq:continous_2} \\
     &\frac{\partial \pi}{\partial t} - \nabla \cdot\Big(\nabla \pi -\frac{2 \nabla H}{H} \pi \Big)=0, \label{eq:continous_3} \\
     &\frac{\partial H}{\partial t} - \frac{1}{\tau} \Big(\rho A e^{- \pi} - H \Big) = 0, \label{eq:continous_4}
\end{align}
where $\eta>0$ characterizes the strength of the
near-repeat victimization and the broken-window effect and $q>0$ is a source term related to criminal regeneration, time scale for repeat victimization, and probability of arrests.
See \cite{Hao2026b} for more details.

\begin{remark}\label{rem:const_param}
    For simplicity, in this paper we assume
    $\Ast$, $\eta$, and $q$ are constant.
    This will allow us to cancel some terms and simplify the presentation of our approaches.
    The work presented here can be easily extended to the case of spatially varying parameters, which we have considered in \cite{Hao2026,Hao2026b}.
\end{remark}

Problem \eqref{eq:continous_1}-\eqref{eq:continous_4} is supplemented with boundary conditions:
\begin{align}
   \nabla A {\cdot} \bn =
   \nabla \rho {\cdot} \bn =
   \nabla \pi {\cdot} \bn &= 0, \quad \text{on } \partial\Omega,\; t>0, \label{eq:BC}
\end{align}
where $\bn$ is the outward unit normal vector.
We prescribe the initial conditions
\begin{align}
     A(0,\bm{x}) =   A_0,\quad
\rho(0,\bm{x}) = \rho_0, \quad
\pi(0,\bm{x}) =
\pi_0, \label{eq:ic}
\end{align}
where
$A_0$, $\rho_0$, and $\pi_0$ are positive constants prescribing the initial level of attractiveness,
the initial criminal density, and the initial police density.
Since $H$ is delayed crime data, which can be estimated by $\rho A e^{-\pi}$, the initial condition for \eqref{eq:continous_4} is given by
\begin{align}
    H(0,\bm{x}) = \rho(0,\bm{x})A(0,\bm{x}) e^{- \pi(0,\bm{x})}.\label{eq:ic4}
\end{align}

We note that both advection terms in \eqref{eq:continous_2} and \eqref{eq:continous_3}
involve logarithmic sensitivities, i.e., $\nabla A / A$ and $\nabla H / H$.
Hence, problem \eqref{eq:continous_1}-\eqref{eq:continous_4} is well defined only as long as
$A$ and $H$ remain positive, while $\rho$ and $\pi$ are expected to remain nonnegative.
Moreover, through the deterrence factor $e^{-\pi}$, the density of police enters
the equations for all other variables.

Since eq.~\eqref{eq:continous_4} is an ODE, no boundary condition is prescribed for $H$. Nevertheless, it is easy to show that $H$ satisfies
\begin{equation}\label{eq:boundary_for_H}
    \nabla H \cdot \bn = 0 \quad \text{on } \partial\Omega,\; t>0,
\end{equation} a property that will be needed in the weak formulation of \eqref{eq:continous_3}.

Taking into account \eqref{eq:boundary_for_H} and boundary conditions \eqref{eq:BC},
the weak problem associated with eq.~\eqref{eq:continous_1}-\eqref{eq:continous_3} reads: Find $A(t,\bx), H(t,\bx)\in W^{1,\infty}(\Omega)$ and $\rho(t,\bx), \pi(t,\bx) \in H^1(\Omega)$ %
such that
\begin{align}
        &\int_\Omega \frac{\partial A}{\partial t} u\;d\Omega + \int_\Omega \eta\nabla A \cdot \nabla u\;d\Omega
        +\int_\Omega Au \;d\Omega
        - \int_\Omega \rho A e^{-\pi} u \;d\Omega  = \int_\Omega \Ast  u \;d\Omega , \label{eq:weak_1} \\
        & \int_\Omega \frac{\partial \rho}{\partial t} v \;d\Omega +\int_\Omega  \Big(\nabla\rho - \frac{2 \nabla A}{A} \rho \Big) \cdot \nabla v\;d\Omega
        +\int_\Omega \rho A e^{-\pi}  v\;d\Omega - \int_\Omega q e^{-\pi} v\;d\Omega=  0,  \label{eq:weak_2} \\
        &\int_\Omega \frac{\partial \pi}{\partial t} w \;d\Omega +\int_\Omega  \Big(\nabla\pi - \frac{2 \nabla H}{H} \pi \Big) \cdot \nabla w\;d\Omega
         = 0, \label{eq:weak_3}
\end{align}
for a.e.~$t\in(0,T)$ and all $(u,v,w)\in H^1(\Omega)^3$, and \eqref{eq:continous_4} holds a.e.~in $(0,T)\times\Omega$.
Note that the term $-\eta \Delta \Ast$ in \eqref{eq:continous_1} vanishes, since $\Ast$ is constant (see Remark~\ref{rem:const_param}).
Because of the terms $\nabla A/A$ and $\nabla H/H$ in \eqref{eq:weak_2} and \eqref{eq:weak_3}, we require $A, H\in W^{1,\infty}(\Omega)$ to be bounded away from zero.

We note that the total police population is conserved ({see~\cite{Hao2026b}}).
Indeed, taking $w=1$ in \eqref{eq:weak_3} gives
\begin{equation}\label{eq:p_cons}
    \int_\Omega \pi(t,\bx)\,d\Omega
    = \int_\Omega \pi(0,\bx)\,d\Omega,
    \qquad t\in[0,T].
\end{equation}

\subsection{Discretization}
We first introduce a time discretization of problem \eqref{eq:continous_4}, \eqref{eq:weak_1}-\eqref{eq:weak_3}. For this,
we set $\Delta t=T/N_t$ for a given positive
integer $N_t$ and $t^n=n\Delta t$ for $n=0,\ldots,N_t$.
We denote by $y^n$ the approximation of a quantity $y$
at time $t^n$.
For simplicity, we apply Backward Differentiation Formula
of order 1 to \eqref{eq:continous_4}, \eqref{eq:weak_1}-\eqref{eq:weak_3} and get the following problem at the generic time $t^{n+1}$: for $n\geq 0$,
given $A^n$, $\rho^n$, $\pi^n$, and $H^n$, find $A^{n+1},H^{n+1}\in W^{1,\infty}(\Omega)$ and $\rho^{n+1}, \pi^{n+1} \in H^1(\Omega) $, %
such that
\begin{align}
           & \frac{1}{\Delta t}\int_\Omega A^{n+1} u\;d\Omega +
           \int_\Omega \eta\nabla A^{n+1} \cdot \nabla u\;d\Omega +\int_\Omega A^{n+1}u \;d\Omega \cl
           & \quad \quad -\int_\Omega \rho^{n+1} A^{n+1} e^{- \pi^{n+1}} u \;d\Omega  = \int_\Omega \Ast u \;d\Omega + \frac{1}{\Delta t}\int_\Omega A^{n} u\;d\Omega, \label{eq:implicit-1} \\
         & \frac{1}{\Delta t} \int_\Omega \rho^{n+1} v \;d\Omega +\int_\Omega  \Big(\nabla\rho^{n+1} - \frac{2\nabla A^{n+1}}{A^{n+1}} \rho^{n+1}\Big) \cdot \nabla v\;d\Omega  \cl
         & \quad \quad +\int_\Omega \rho^{n+1} A^{n+1} e^{-\pi^{n+1}}  v\;d\Omega - \int_\Omega q e^{-\pi^{n+1}} v\;d\Omega = \frac{1}{\Delta t} \int_\Omega \rho^{n} v \;d\Omega, \label{eq:implicit-2} \\
         & \frac{1}{\Delta t} \int_\Omega \pi^{n+1} w \;d\Omega +\int_\Omega  \Big(\nabla\pi^{n+1} - \frac{2\nabla H^{n+1}}{H^{n+1}} \pi^{n+1}\Big) \cdot \nabla w\;d\Omega =  \frac{1}{\Delta t} \int_\Omega \pi^{n} w \;d\Omega, \label{eq:implicit-3} \\
         & \frac{1}{\Delta t} H^{n+1}   -  \frac{1}{\tau} \Big(\rho^{n+1} A^{n+1} e^{- \pi^{n+1}} - H^{n+1} \Big)   = \frac{1}{\Delta t}  H^{n}, \label{eq:implicit-4}
\end{align}
for all $(u,v,w)\in H^1(\Omega)^3$.
System \eqref{eq:implicit-1}-\eqref{eq:implicit-4} is
implicit and coupled. Once discretized in space, it yields a nonlinear system whose size grows as the mesh is refined
and which needs to be solved at each time step through an iterative method such as Newton’s method. Since
problem \eqref{eq:implicit-1}-\eqref{eq:implicit-4} features several strongly nonlinear terms, the convergence of such an
iterative method may be slow, especially for large-scale simulations \cite{quarteroni1994numerical,saad2003iterative}.

Below, we present two strategies that avoid solving the entire nonlinear system at each time step. The first is the fully decoupled approach used in~\cite{Hao2026b}, which solves for one variable at a time. The second is a coupled approach that solves for $A$ and $\rho$ simultaneously, with the reaction term $\rho A$ linearized by Newton's method. To speed up both iterations, we combine them with Anderson acceleration, which we introduce in Sec.~\ref{sec:AA}. Before presenting both approaches, let us introduce details about the spatial discretization.

Let $\mathcal{T}_h$ be a partition of $\Omega$ into $N_e$
quadrilaterals $Q_k$ such that
\begin{equation*}
    \overline{\Omega} = \bigcup_{k = 1}^{N_e} \overline{Q}_{k},
\end{equation*}
with
\begin{equation*}
    h = \max_{Q_k \in \mathcal{T}_h} h_k, \quad h_k = \text{diam}(Q_k), \quad k= 1, \dots, N_e.
\end{equation*}
Let $\mathbb{Q}_1 (Q_k)$ be the set of polynomials, defined on $Q_k$, of degree less than or equal to one in each space variable, and let $V_h (\Omega)$ be the space of globally continuous functions on $\overline{\Omega}$
whose restriction to each $Q_k \in \mathcal{T}_h$ belongs to $\mathbb{Q}_1 (Q_k)$:
\begin{equation*}
    V_h (\Omega) = \{ v_h \in C^0 (\overline{\Omega})~:~v_h |_{Q_k} \in \mathbb{Q}_1(Q_k), \forall Q_k \in \mathcal{T}_h  \}.
\end{equation*}
We denote by $\{ \phi_i \}_{i = 1}^{N_Q}$ a basis for $V_h$, where $N_Q$ is the total number of degrees of freedom in $\Omega$. Note that property \eqref{eq:p_cons} is preserved by the
Backward Euler-Galerkin discretization, since
constant functions belong to $V_h(\Omega)$.

Finally, we introduce some matrices that will be used in Secs.~\ref{sec:decoupled} and
\ref{sec:coupled}.
Let $M$ be the mass matrix, $K$ the stiffness matrix, $N(\omega)$ a weighted mass matrix
and $D(\omega)$ a ``weighted'' divergence matrix, for a given bounded weight function $\omega$, with $\omega > 0$ in $\overline{\Omega}$ and $\nabla \omega$ bounded in the case of $D(\omega)$:
\begin{align}
    &M_{i,j}= \int_\Omega \phi_j \phi_i \;d\Omega, \quad
    K_{i,j} = \int_\Omega \nabla \phi_j \cdot \nabla\phi_i \;d\Omega, \label{eq:M_K} \\
    &N(\omega)_{i,j} = \int_\Omega \omega \phi_j \phi_i \;d\Omega, \quad
    D(\omega)_{i,j} = \int_\Omega \phi_j \, \frac{2 \nabla \omega}{\omega} \cdot \nabla \phi_i \;d\Omega, \label{eq:N_D}
\end{align}
for $i, j = 1, \dots, N_Q$. All the integrals are computed with a quadrature rule of sufficiently high order.  In particular, when the weight $\omega$ is a nonlinear function of finite element functions,
e.g., $\omega = \rho_h e^{-\pi_h}$, it is evaluated at the quadrature points.

\section{Fixed-point iterations and Anderson acceleration}\label{sec:AA}

Both approaches considered in this paper
solve the space discretized  nonlinear problem \eqref{eq:implicit-1}-\eqref{eq:implicit-4} at each time step iteratively. As we will show in Sec.~\ref{sec:AA-decoupled} and \ref{sec:AA-coupled}, both iterations can be written as
\begin{equation}\label{eq:fixed-point}
    \mathbf{u}^{k+1} = G(\mathbf{u}^{k}), \quad k \geq 0,
\end{equation}
for a suitable map $G: \mathbb{R}^{\ell N_Q} \to \mathbb{R}^{\ell N_Q}$, where vector $\mathbf{u}$ collects the nodal values of $\ell$ model variables and each evaluation of $G$
requires the solution of one or more linear systems. We will denote by $G_d$ and $G_c$ the maps associated with the fully decoupled approach in Sec.~\ref{sec:AA-decoupled} and the coupled approach in Sec.~\ref{sec:AA-coupled}, respectively.

\begin{remark}\label{rem:notation}
{Throughout the paper, the superscript $n$ denotes the
time level, whereas the superscript $k$ denotes the
fixed-point iteration index within a given time step.
For simplicity, we omit the time-level index from
the iterates, since all fixed-point iterations are
performed at a fixed time step.}

\end{remark}

Fixed-point iterations like \eqref{eq:fixed-point} converge linearly, with asymptotic rate
given by the spectral radius of the Jacobian of $G$ at the fixed point, which is independent
of the norm.
In practice, this rate can be estimated from the fixed-point residual $\mathbf{f}^{k} = G(\mathbf{u}^{k}) - \mathbf{u}^{k}$.
Since accelerated iterations, such as those introduced below, converge R-linearly \cite{TothKelley2015}, i.e., the residual is bounded by a geometrically decreasing sequence
but need not decrease at every iteration, we do not use the ratio of two consecutive residuals. Instead, we use the average contraction factor
\begin{equation}\label{eq:contraction-factor}
    \bar\mu = \left( \frac{\|\mathbf{f}^{k_1}\|}{\|\mathbf{f}^{k_0}\|} \right)^{1/(k_1 - k_0)}, \quad 0 \le k_0 < k_1,
\end{equation}
i.e., the geometric mean of the reduction of the residual between iterations $k_0$ and $k_1$, where $\|\cdot\|$ can be any norm in $\mathbb{R}^{\ell N_Q}$.
This quantity estimates the R-linear convergence factor, and it is meaningful for both plain and accelerated iterations.
For the plain iteration \eqref{eq:fixed-point}, $\mathbf{f}^{k} = \mathbf{u}^{k+1} - \mathbf{u}^{k}$, and $\bar\mu$ typically approaches the spectral radius of the Jacobian as $k_0$ and $k_1$ increase.
Convergence is slow when this spectral radius is close to one, and the iteration
may fail to converge when it exceeds one. For the problem at hand, the spectral radius depends
on the time step, on the model parameters, on the current state (through the Jacobian at the fixed point), and on how each approach handles the coupling among the model variables,
as we will show in Sec.~\ref{sec:contraction} and \ref{sec:contraction-coupled}.

Anderson acceleration \cite{Anderson1965,Walker2011,he2026generalizedAA,he2026modifiedAA} is designed to speed up iterations like \eqref{eq:fixed-point}.
Recall that $\mathbf{f}^{k} = G(\mathbf{u}^{k}) - \mathbf{u}^{k}$ denotes the fixed-point residual.
Instead of taking $\mathbf{u}^{k+1} = G(\mathbf{u}^{k})$, AA with depth $m \geq 0$, denoted AA($m$), uses the last $m_k + 1$ evaluations of $G$, with $m_k = \min\{k, m\}$:
\begin{equation}\label{eq:AA-update}
    \mathbf{u}^{k+1} = G(\mathbf{u}^{k}) - \mathcal{G}^k \boldsymbol{\gamma}^{k}, \quad
    \boldsymbol{\gamma}^{k} = \operatorname*{arg\,min}_{\boldsymbol{\gamma} \in \mathbb{R}^{m_k}} \left\| \mathbf{f}^{k} - \mathcal{F}^k \boldsymbol{\gamma} \right\|_B,
\end{equation}
where
\begin{align*}
    \mathcal{G}^k & = \left[ \Delta G^{k-m_k}, \dots, \Delta G^{k-1} \right], \quad \Delta G^i  = G(\mathbf{u}^{i+1}) - G(\mathbf{u}^{i}), \\
    \mathcal{F}^k & = \left[ \Delta \mathbf{f}^{k-m_k}, \dots, \Delta \mathbf{f}^{k-1} \right], \qquad \Delta \mathbf{f}^{i}  = \mathbf{f}^{i+1} - \mathbf{f}^{i}.
\end{align*}
Here, $\|\mathbf{x}\|_B = (\mathbf{x}^T B \mathbf{x})^{1/2}$ for a given symmetric positive definite matrix $B$. The standard choice is $B = I$, i.e., the Euclidean norm.
We discuss the choice of $B$ below. Note that $\boldsymbol{\gamma}^{k}$ minimizes the norm of a linear combination of the last $m_k + 1$ residuals, and the same combination is applied to the corresponding evaluations of $G$.
For $m = 0$, \eqref{eq:AA-update} reduces to the original iteration \eqref{eq:fixed-point}.

The problem for $\boldsymbol{\gamma}^{k}$ in \eqref{eq:AA-update} is a linear least-squares problem.
Since $\mathcal{F}^k \in \mathbb{R}^{\ell N_Q \times m_k}$ with {$m_k \ll \ell N_Q$}, it is highly overdetermined. %
If $\mathcal{F}^k$ has full column rank, the solution is unique and satisfies the normal equations
\begin{equation*}
    (\mathcal{F}^k)^T B \, \mathcal{F}^k \, \boldsymbol{\gamma}^{k} = (\mathcal{F}^k)^T B \, \mathbf{f}^{k}.
\end{equation*}
{In our implementation, we form the $m_k \times m_k$ Gram matrix $(\mathcal{F}^k)^T B \, \mathcal{F}^k$ and solve the normal equations directly. This is inexpensive since $m_k$ is small, and, if the products of $B$ with the columns of $\mathcal{F}^k$ are stored, it requires only a few sparse matrix-vector products with $B$ per iteration.
Since consecutive differences $\Delta \mathbf{f}^{i}$ can become nearly linearly dependent as the iteration converges, $\mathcal{F}^k$ may be ill-conditioned, and forming the Gram matrix squares its condition number.
Therefore, we drop the oldest columns of $\mathcal{F}^k$ and $\mathcal{G}^k$ as soon as the condition number of the Gram matrix exceeds a threshold ($10^{10}$), similarly to the column deletion in \cite{Walker2011}.
A more stable, but more involved, alternative is to compute a thin QR factorization of $\mathcal{F}^k$ with respect to the inner product induced by $B$ and solve a small triangular system \cite{Walker2011}.
{Other techniques for solving the least-squares problem can be found, e.g., in \cite{pollock2023filtering,lockhart2022performance}.}}

Alg.~\ref{alg:AA-generic} summarizes one time step of a fixed-point iteration \eqref{eq:fixed-point} accelerated with AA($m$).
The evaluation of $G$ in line 2 is specified in Sec.~\ref{sec:AA-decoupled} for $G = G_d$ and in Sec.~\ref{sec:AA-coupled} for $G = G_c$.
Alg.~\ref{alg:AA-generic} uses a stopping criterion based on the relative fixed-point residual, i.e., the relative difference between $\mathbf{u}^{k}$ and $G(\mathbf{u}^{k})$, measured in a suitable norm $\|\cdot\|$.
The criterion is checked before the AA update, so that, once a plain fixed-point step satisfies it, no further acceleration is performed and the accepted solution is always an evaluation of $G$.
The specific criterion used for our problem is given in Sec.~\ref{sec:decoupled}.
In Alg.~\ref{alg:AA-generic}, the iteration is initialized with the solution at the previous time level, $\mathbf{u}^{0} = \mathbf{u}^{n}$.
A linear extrapolation from the two previous time levels, $\mathbf{u}^{0} = 2\mathbf{u}^{n} - \mathbf{u}^{n-1}$, may reduce the initial error for smooth solutions.

\begin{algorithm}
\caption{{AA($m$)-accelerated fixed-point iteration \eqref{eq:fixed-point}}: time step from $t^n$ to $t^{n+1}$}\label{alg:AA-generic}
\begin{algorithmic}[1]
\Require set $\mathbf{u}^{0} = \mathbf{u}^{n}$; depth $m \geq 0$; tolerance $\texttt{tol}$; maximum number of iterations $k_{\max}$
\For{$k = 0, 1, \dots, k_{\max} - 1$}
    \State Evaluate $G(\mathbf{u}^{k})$ and compute $\mathbf{f}^{k} = G(\mathbf{u}^{k}) - \mathbf{u}^{k}$
    \If{$\|\mathbf{f}^{k}\| \le \texttt{tol} \, \|\mathbf{u}^{k}\|$}
        \State Set $\mathbf{u}^{n+1} = G(\mathbf{u}^{k})$ and \textbf{stop}
        \Comment{converged with a plain fixed-point step}
    \EndIf
    \State Set $m_k = \min\{k, m\}$
    \If{$m_k = 0$}
        \State $\mathbf{u}^{k+1} = G(\mathbf{u}^{k})$
    \Else
        \State $\mathcal{F}^k = [\Delta \mathbf{f}^{k-m_k}, \dots, \Delta \mathbf{f}^{k-1}]$, \ $\mathcal{G}^k = [\Delta G^{k-m_k}, \dots, \Delta G^{k-1}]$
        \State $\boldsymbol{\gamma}^{k} = \operatorname*{arg\,min}_{\boldsymbol{\gamma} \in \mathbb{R}^{m_k}} \| \mathbf{f}^{k} - \mathcal{F}^k \boldsymbol{\gamma} \|_B$
        \State $\mathbf{u}^{k+1} = G(\mathbf{u}^{k}) - \mathcal{G}^k \boldsymbol{\gamma}^{k}$
    \EndIf
\EndFor
\State \textbf{if} $k_{\max}$ is reached, the time step is considered failed
\end{algorithmic}
\end{algorithm}

Each iteration of Alg.~\ref{alg:AA-generic} requires one evaluation of $G$, i.e., the same linear solves as one iteration of the original fixed-point iteration.
The extra work consists of storing the last $m_k$ differences $\Delta \mathbf{f}^{i}$ and $\Delta G^i$ and solving a least-squares problem with only {$m_k \ll \ell N_Q$} unknowns,
which is negligible compared to the linear solves.
Moreover, AA treats $G$ as a black box, so the linear solvers used in the evaluation of $G$ are unaffected.

The unknown $\mathbf{u}$, and hence $\mathbf{f}^{k}$, consists of blocks associated with different model variables.
With $B = I$, if the variables have very different magnitudes, the larger blocks dominate the minimization in \eqref{eq:AA-update} and the acceleration
acts mostly on the corresponding variables. This can be avoided by choosing a block-diagonal matrix $B$ that weights the blocks.
More generally, the least-squares problem may be posed in any norm without affecting the local convergence of AA \cite{TothKelley2015}, and a suitable choice of norm can improve its performance (see, e.g., \cite{hawkins2026choice,YangTownsendAppelo2022}). In \cite{YangTownsendAppelo2022} a Sobolev norm is used, which becomes a weighted norm after discretization.
Moreover, $B$ can be chosen so that $\|\cdot\|_B$ is the norm in which $G$ is contractive, as in \cite{Pollock2019,Evans2020}, so that the results in \cite{Evans2020} apply.
{In this setting, $G$ maps $X = V_h(\Omega)^\ell$, identified with $\mathbb{R}^{\ell N_Q}$, into itself,
and the norm $\|\cdot\|_B$ is used both for the contraction of $G$ and for the residual, as in \cite{Evans2020}.}
In our computations, we choose
\begin{equation}\label{eq:B}
    B = \text{diag} \left( w_1^2 (M + \Delta t K), \dots, w_\ell^2 (M + \Delta t K) \right),
\end{equation}
where $w_j = 1 / \max \{ \|\mathbf{v}_j^n\|_2, \epsilon \}$, with $\mathbf{v}_j^n$ the vector of nodal values of the $j$-th variable in $\mathbf{u}$ at time $t^n$
and $\epsilon > 0$ a small constant that prevents division by zero. The weights are fixed for the whole time step.
With this choice, $\|\cdot\|_B$ is a block-weighted version of the norm in which we show that $G_d$ is a contraction (see Sec.~\ref{sec:contraction}),
so that the results in \cite{Evans2020} apply to $G_d$ (see Remark~\ref{rem:AA-norm}). Instead, in general $G_c$ is  not a contraction in this norm (see Sec.~\ref{sec:contraction-coupled}).

The convergence of AA has been analyzed mostly for contractive maps.
If $G$ is a contraction that is continuously differentiable with Lipschitz continuous derivative in a neighborhood of the fixed point,
AA($m$) converges locally R-linearly, { provided that the coefficients in the Anderson combination remain bounded}~\cite{TothKelley2015}. Moreover, AA($m$) improves the linear convergence rate of \eqref{eq:fixed-point} by the gain factor of the least-squares problem,
up to higher-order terms \cite{Evans2020}, provided that the least-squares problem is posed in the inner-product norm in which $G$ is contractive.
Results for noncontractive maps can be found in \cite{pollock2021anderson}.  One of the key advantages of AA is its ability to promote convergence even for noncontractive fixed-point iterations.
Asymptotically, AA behaves like GMRES applied to the problem linearized around the fixed point
\cite{Walker2011,DeSterckHe2021}. Hence, while the asymptotic rate of the plain iteration
\eqref{eq:fixed-point} is determined by the spectral radius of the Jacobian of $G$ at the fixed point,
the convergence of AA depends on the whole spectrum of the Jacobian, e.g., on the clustering of its eigenvalues.
In particular, AA can converge fast even when the spectral radius is close to one.

\section{A fully decoupled approach}\label{sec:decoupled}

In \cite{Hao2026b}, we adopted an approach that
decouples system \eqref{eq:implicit-1}-\eqref{eq:implicit-4} while linearizing it
and performs fixed-point iterations on linear subproblems to avoid possible instabilities for certain combinations of parameter values.
We directly state it in its fully discretized form, with the finite element space $V_h(\Omega)$ introduced in Sec.~\ref{sec:pd}.
The algorithm is as follows:
at every time step $t^{n+1}$, given $A^n_h$, $\rho^n_h$, $\pi^n_h$, and $H^n_h$, set $A_h^0 = A_h^n$,
$\rho^{0}_h = \rho^{n}_h$ and $\pi^{0}_h = \pi^{n}_h$ (see Sec.~\ref{sec:AA}) and
perform at iteration $k+1$, $k \ge 0$ (see notation explanation in Remark~\ref{rem:notation}):
\begin{itemize}
    \item[-] \emph{Step 1}:  find $A^{k+1}_h \in V_h (\Omega)$ such that
    \begin{align}
           &\frac{1}{\Delta t}\int_\Omega A^{k+1}_h u_h\;d\Omega +
           \int_\Omega \eta\nabla A^{k+1}_h \cdot \nabla u_h\;d\Omega +\int_\Omega A^{k+1}_h u_h \;d\Omega \cl
           & \quad \quad -\int_\Omega \rho^{k}_h {e^{- \pi^{k}_h}} A^{k+1}_h u_h \;d\Omega  = \int_\Omega \Ast u_h \;d\Omega  + \frac{1}{\Delta t}\int_\Omega A^{n}_h u_h\;d\Omega, \label{eq:iterative-1}
\end{align}
for all $u_h \in V_h (\Omega)$.
 \item[-] \emph{Step 2}:  find $\rho^{k+1}_h\in V_h (\Omega)$ such that
 \begin{align}
          & \frac{1}{\Delta t} \int_\Omega \rho^{k+1}_h v_h \;d\Omega +\int_\Omega  \Big(\nabla\rho^{k+1}_h - \frac{2\nabla A^{k+1}_h}{A^{k+1}_h} \rho^{k+1}_h\Big) \cdot \nabla v_h\;d\Omega  \cl
         & \quad \quad +\int_\Omega \rho^{k+1}_h A^{k+1}_h {e^{-\pi^{k}_h}} v_h\;d\Omega = \int_\Omega q {e^{-\pi^{k}_h}} v_h\;d\Omega + \frac{1}{\Delta t} \int_\Omega \rho^{n}_h v_h \;d\Omega, \label{eq:iterative-2}
\end{align}
for all $v_h \in V_h (\Omega)$.
\item[-] \emph{Step 3}:  set $H^{k+1}_h \in V_h(\Omega)$ such that
\begin{align}
 \left( \frac{1}{\Delta t} + \frac{1}{\tau} \right)  H^{k+1}_h    =  \frac{1}{\tau} \Big(\rho^{k+1}_h A^{k+1}_h e^{- \pi^{k}_h} \Big)  +  \frac{1}{\Delta t}  H^{n}_h \label{eq:iterative-3}
\end{align}
holds at the nodes of $\mathcal{T}_h$.
\item[-] \emph{Step 4}:  find $\pi^{k+1}_h\in V_h(\Omega)$ such that
 \begin{align}
         & \frac{1}{\Delta t} \int_\Omega \pi^{k+1}_h w_h \;d\Omega +\int_\Omega  \Big(\nabla\pi^{k+1}_h - \frac{2\nabla H^{k+1}_h}{H^{k+1}_h} \pi^{k+1}_h\Big) \cdot \nabla w_h\;d\Omega  =  \frac{1}{\Delta t} \int_\Omega \pi^{n}_h w_h \;d\Omega, \label{eq:iterative-4}
\end{align}
for all $w_h \in V_h (\Omega)$.
\item[-] \emph{Step 5}: Check the stopping criterion, e.g.
\begin{equation}\label{eq:crit}
    \frac{\|A^{k+1}_h - A^{k}_h\|}{\|A^{k}_h\|}< \texttt{tol}_1, \quad  \frac{\|\rho^{k+1}_h -\rho^{k}_h \|}{\|\rho^{k}_h \|}< \texttt{tol}_2,
    \quad \frac{\|\pi^{k+1}_h-\pi^{k}_h\|}{\|\pi^{k}_h\|}< \texttt{tol}_3,
\end{equation}
where $\|\cdot\|$ denotes the norm of $L^2(\Omega)$ and $\texttt{tol}_1$, $\texttt{tol}_2$, and $\texttt{tol}_3$ are given stopping tolerances. If the stopping criterion is not satisfied, repeat steps 1-5. If it is satisfied, set $A^{n+1}_h = A^{k+1}_h$, $\rho^{n+1}_h = \rho^{k+1}_h$, $\pi^{n+1}_h = \pi^{k+1}_h$,
$H^{n+1}_h = H^{k+1}_h$.
\end{itemize}

At a given time step, algorithm \eqref{eq:iterative-1}-\eqref{eq:crit} iterates
over four subproblems until criterion \eqref{eq:crit} is satisfied: three linear
problems, \eqref{eq:iterative-1}, \eqref{eq:iterative-2}, and \eqref{eq:iterative-4},
and the explicit update \eqref{eq:iterative-3}. Since the equation for $H$
contains no spatial derivatives, \eqref{eq:iterative-3} is a nodewise update and
requires no linear solve.
In practice, at step 5 we use the analogous criterion for the vectors of nodal values, measured in the Euclidean norm.
Since this criterion is based on relative differences, on quasi-uniform meshes it is equivalent to \eqref{eq:crit} up to constants that depend only on the mesh quality.

Next, we write \eqref{eq:iterative-1}-\eqref{eq:iterative-2},\eqref{eq:iterative-4} in matrix form.
For this, we use the matrices defined in \eqref{eq:M_K} and \eqref{eq:N_D}. Moreover,
let $\bA^{k+1}, \br^{k+1}, \bp^{k+1}, \bH^{k+1} \in \mathbb{R}^{N_Q}$ be the vectors of
the nodal values of $A^{k+1}_h$, $\rho^{k+1}_h$,
$\pi^{k+1}_h$, and $H^{k+1}_h$.
In the matrix forms below, the nonlinear coefficients, e.g., $\rho^{k}_h e^{-\pi^{k}_h}$, enter the matrices $N(\cdot)$ and $D(\cdot)$ defined in \eqref{eq:N_D} as functions,
i.e., they are evaluated at the quadrature points. The only exception is the update of $H$, which is performed with respect to the nodal values.

The matrix form of \eqref{eq:iterative-1}
is given by: find the vector $\bA^{k+1} \in \mathbb{R}^{N_Q}$ such that
\begin{equation}\label{eq:mat1bis}
     \left[ \frac{1 + \Delta t}{\Delta t} M + \eta K - N(c^{k}_h) \right] \bA^{k+1} = \bb_A^{k+1},
\end{equation}
where $c^{k}_h = \rho^{k}_h e^{- \pi^{k}_h}$ and
the right-hand side vector is given by
\begin{equation*}
    \bb_A^{k+1} =  M \bA^{\mathrm{st}} + \frac{1}{\Delta t} M \bA^{n},
\end{equation*}
with $\bA^{\mathrm{st}} \in \mathbb{R}^{N_Q}$ the vector of the nodal values of $\Ast$.

The matrix form of \eqref{eq:iterative-2}
is given by: find the vector $\br^{k+1} \in \mathbb{R}^{N_Q}$ such that
\begin{equation}\label{eq:mat2bis}
    \left[ \frac{1}{\Delta t} M + K - D(A^{k+1}_h) + N(A^{k+1}_h e^{-\pi^{k}_h}) \right] \br^{k+1} = \bb_\rho^{k+1},
\end{equation}
where
the vector at the right-hand side $\bb_\rho^{k+1}$ has $i$-th component
\begin{equation*}
    \bb_{\rho,i}^{k+1} = \int_\Omega q e^{-\pi^{k}_h} \phi_i \;d\Omega + \frac{1}{\Delta t} (M\br^{n})_i,
\end{equation*}
for $i = 1, \dots, N_Q$. We note that
$A^{k+1}_h$ is known from solving \eqref{eq:mat1bis}.

In terms of the vectors of nodal values, \eqref{eq:iterative-3} reads: set $\bH^{k+1} \in \mathbb{R}^{N_Q}$ as follows
\begin{align}
  \bH^{k+1}    =  \frac{\Delta t}{\tau + \Delta t} \bd^{k+1}  +  \frac{\tau}{\tau + \Delta t}  \bH^{n}, \label{eq:iterative-3-v}
\end{align}
where $\bd^{k+1}$ is the vector of the nodal values of function $\rho^{k+1}_h A^{k+1}_h e^{- \pi^{k}_h}$.

Finally, the matrix form of \eqref{eq:iterative-4}
is given by: find the vector $\bp^{k+1} \in \mathbb{R}^{N_Q}$ such that
\begin{equation}\label{eq:mat4}
    \left[ \frac{1}{\Delta t} M + K - D(H^{k+1}_h) \right] \bp^{k+1} = \frac{1}{\Delta t} M \bp^{n},
\end{equation}
where $H^{k+1}_h$ is known from \eqref{eq:iterative-3-v}.

The three linear systems \eqref{eq:mat1bis}, \eqref{eq:mat2bis}, and \eqref{eq:mat4} have different algebraic properties.
Since $M$, $K$, and $N(c^{k}_h)$ are symmetric, the matrix in \eqref{eq:mat1bis} is symmetric. Moreover, since $\bx^T N(c^{k}_h) \bx \le \| c^k_h \|_{L^\infty(\Omega)} \, \bx^T M \bx$ for all $\bx \in \mathbb{R}^{N_Q}$
and $K$ is symmetric positive semidefinite (because of the homogeneous Neumann boundary conditions, its kernel consists of the constant vectors),
a sufficient condition for this matrix to be symmetric positive definite is
\begin{equation}\label{eq:suff_cond}
    1 + \frac{1}{\Delta t} > \| c^k_h \|_{L^\infty(\Omega)},
\end{equation}
i.e., the time step has to be small enough with respect to the size of $\rho^k_h e^{-\pi^k_h}$.
On the other hand, the matrices in \eqref{eq:mat2bis} and \eqref{eq:mat4} are nonsymmetric because of the weighted divergence matrices $D(A^{k+1}_h)$ and $D(H^{k+1}_h)$.
In fact, \eqref{eq:mat2bis} and \eqref{eq:mat4} are discretizations of advection-diffusion-reaction problems, with advection fields $2\nabla A^{k+1}_h/A^{k+1}_h$ and $2\nabla H^{k+1}_h/H^{k+1}_h$, respectively.
Therefore, they call for Krylov methods for nonsymmetric systems, such as GMRES or BiCGStab \cite{saad2003iterative}, and become more challenging to solve in regions where $|\nabla A^{k+1}_h|/A^{k+1}_h$ or $|\nabla H^{k+1}_h|/H^{k+1}_h$ is large, i.e., where advection dominates diffusion.
See Sec.~\ref{sec:num_res} for more details.

\subsection{Fixed-point formulation and Anderson acceleration}
\label{sec:AA-decoupled}

Algorithm \eqref{eq:iterative-1}-\eqref{eq:crit} is a nonlinear block Gauss--Seidel iteration:
at each iteration, the subproblems are solved in sequence, each using the most recent approximation of the other variables.
In terms of the vectors of nodal values, i.e., in the form \eqref{eq:mat1bis}, \eqref{eq:mat2bis}, \eqref{eq:iterative-3-v}, \eqref{eq:mat4},
it can be written as a fixed-point iteration.  Given the solution at time $t^n$,
$\bA^{k+1}$ depends only on $\br^{k}$ and $\bp^{k}$, through $N(c^{k}_h)$, with $c^{k}_h = \rho^{k}_h e^{-\pi^{k}_h}$. Vector $\br^{k+1}$ depends only on $\bA^{k+1}$ and $\bp^{k}$. Vector
$\bH^{k+1}$ depends only on $\bA^{k+1}$, $\br^{k+1}$, and $\bp^{k}$ through vector $\bd^{k+1}$. Finally, $\bp^{k+1}$ depends only on $\bH^{k+1}$.
Hence, iteration $k+1$ is fully determined by $\br^{k}$ and $\bp^{k}$, while $\bA^{k}$ and $\bH^{k}$ play no role in it.

Let $\mathbf{u} = (\br, \bp) \in \mathbb{R}^{2N_Q}$ and let $G_d: \mathbb{R}^{2N_Q} \to \mathbb{R}^{2N_Q}$ be the map
that takes $\mathbf{u}^{k} = (\br^{k}, \bp^{k})$ to $\mathbf{u}^{k+1} = (\br^{k+1}, \bp^{k+1})$ through steps 1-4 described above.
Then, the fully decoupled algorithm takes the form \eqref{eq:fixed-point} with $G = G_d$ and $\mathbf{u}^{0} = (\br^{n}, \bp^{n})$ (see notation explanation in Remark~\ref{rem:notation}), and its limit is a fixed point of $G_d$.
Each evaluation of $G_d$ consists of the following operations:
\begin{enumerate}
    \item assemble $N(c^{k}_h)$, with $c^{k}_h = \rho^{k}_h e^{-\pi^{k}_h}$, and solve \eqref{eq:mat1bis} for $\bA^{k+1}$;
    \item assemble $D(A^{k+1}_h)$ and $N(A^{k+1}_h e^{-\pi^{k}_h})$, and solve \eqref{eq:mat2bis} for $\widetilde{\br}^{k+1}$;
    \item compute $\bd^{k+1}$ from $\bA^{k+1}$, $\widetilde{\br}^{k+1}$, and $\bp^{k}$, and compute $\bH^{k+1}$ with \eqref{eq:iterative-3-v};
    \item solve \eqref{eq:mat4} for $\widetilde{\bp}^{k+1}$, and set $G_d(\mathbf{u}^{k}) = (\widetilde{\br}^{k+1}, \widetilde{\bp}^{k+1})$.
\end{enumerate}
The tilde denotes the output of the linear solves, which differs from the new iterate $\mathbf{u}^{k+1}$ when AA is applied.
The fully decoupled algorithm with AA($m$) is given by Alg.~\ref{alg:AA-generic} with $G = G_d$, where the stopping criterion is replaced by \eqref{eq:crit} applied to $G_d(\mathbf{u}^{k})$ and $\mathbf{u}^{k}$.
Since $\bA$ and $\bH$ are not part of $\mathbf{u}$, upon convergence $\bA^{n+1}$ and $\bH^{n+1}$ are the vectors computed in the last evaluation of $G_d$.
Each evaluation of $G_d$ requires the solution of the three linear systems \eqref{eq:mat1bis}, \eqref{eq:mat2bis}, and \eqref{eq:mat4}, plus the nodal update \eqref{eq:iterative-3-v}.

\begin{remark}\label{rem:AA-positivity}
The attractiveness $A$ and the crime data $H$ are positive.
The weighted divergence matrices in \eqref{eq:mat2bis} and \eqref{eq:mat4} also require $A^{k+1}_h > 0$ and $H^{k+1}_h > 0$.
Likewise, $\rho$ and $\pi$ should be nonnegative.
Positivity of the discrete solution is not guaranteed even without acceleration: as we will see in Sec.~\ref{sec:num_bounds},
$\pi_h$ takes slightly negative values in some cases, because the Galerkin discretization of the $\pi$ equation does not satisfy
a discrete maximum principle when  advection dominates. In addition, since the coefficients $\boldsymbol{\gamma}^{k}$ can be
negative or large, the accelerated iterate $\mathbf{u}^{k+1}$ in \eqref{eq:AA-update} is not a convex combination of previous
evaluations of $G_d$, and AA may introduce further violations of positivity.
A possible safeguard is to check the sign of the nodal values of $\mathbf{u}^{k+1}$ and, if the check fails, discard the accelerated iterate,
set $\mathbf{u}^{k+1} = G_d(\mathbf{u}^{k})$, and restart AA by clearing $\mathcal{F}^k$ and $\mathcal{G}^k$.
\end{remark}

\subsection{Convergence of the fully decoupled iteration}\label{sec:contraction}

The convergence results for AA recalled in Sec.~\ref{sec:AA} assume that the underlying fixed-point map is contractive
and continuously differentiable with Lipschitz continuous derivative \cite{TothKelley2015,Evans2020}.
In this section, we show that the map $G_d$ defined in Sec.~\ref{sec:AA-decoupled} satisfies these assumptions for a sufficiently small time step (in Proposition \ref{prop:contraction}),
provided that the iterates satisfy suitable bounds. Since proving such bounds for the discrete problem is beyond the scope of this paper,
we state them as an assumption (i.e., Assumption \ref{assump:bounds}) and verify them numerically in Sec.~\ref{sec:num_bounds}.

We identify $\mathbf{u} = (\br, \bp)$ with the pair of finite element functions $(\rho_h, \pi_h) \in X = V_h(\Omega)^2$.
In the analysis, we use the time-step-weighted norm
\begin{equation}\label{eq:tnorm}
    \tnorm{v_h}^2 = \|v_h\|^2 + \Delta t \, \|\nabla v_h\|^2, \qquad
    \tnorm{(\rho_h, \pi_h)}^2 = \tnorm{\rho_h}^2 + \tnorm{\pi_h}^2,
\end{equation}
which is natural for the backward Euler scheme. For a fixed mesh and time step, \eqref{eq:tnorm} is equivalent to the norm used in stopping criterion \eqref{eq:crit},
and the asymptotic convergence rate of a fixed-point iteration does not depend on the choice of norm. Moreover, \eqref{eq:tnorm} is induced by the inner product associated with
the matrix $M + \Delta t K$, which makes it compatible with the least-squares problem in AA (see Remark~\ref{rem:AA-norm} below).
Given two inputs
$\mathbf{u}$ and $\tilde{\mathbf{u}}$ of a map, we denote by a tilde the quantities computed
from $\tilde{\mathbf{u}}$ and with $\delta$ the difference between corresponding quantities,
e.g., $\delta A_h = A_h - \tilde{A}_h$.

All the linear solves in the algorithms in Secs.~\ref{sec:decoupled} and \ref{sec:coupled} involve advection-diffusion-reaction problems with a mass term scaled by $1/\Delta t$.
The following elementary estimate is the main tool of our analysis.

\begin{lem}\label{lem:energy}
Let $\kappa_0 > 0$, $\boldsymbol{\beta} \in L^\infty(\Omega)^2$ with $|\boldsymbol{\beta}| \le b$ in $\Omega$, and $\sigma \in L^\infty(\Omega)$ with $\sigma \ge -s$ in $\Omega$, for some $b, s \ge 0$.
Let $f \in L^2(\Omega)$, $\mathbf{g} \in L^2(\Omega)^2$, and let $w_h \in V_h(\Omega)$ satisfy
\begin{equation}\label{eq:model-problem}
    \frac{1}{\Delta t} (w_h, v_h) + \kappa_0 (\nabla w_h, \nabla v_h) - (\boldsymbol{\beta} w_h, \nabla v_h) + (\sigma w_h, v_h) = (f, v_h) + (\mathbf{g}, \nabla v_h)
\end{equation}
for all $v_h \in V_h(\Omega)$,
where $(\cdot,\cdot)$ denotes the $L^2(\Omega)$ inner product. If $\Delta t \le (b^2 / \kappa_0 + 2 s)^{-1}$, then
\begin{equation}\label{eq:energy}
    \tnorm{w_h} \le \frac{2}{\min\{1, \kappa_0\}} \left( \Delta t \, \|f\| + \Delta t^{1/2} \, \|\mathbf{g}\| \right).
\end{equation}
\end{lem}

\begin{proof}
Take $v_h = w_h$ in \eqref{eq:model-problem}. By Young's inequality, $|(\boldsymbol{\beta} w_h, \nabla w_h)| \le \frac{\kappa_0}{2} \|\nabla w_h\|^2 + \frac{b^2}{2 \kappa_0} \|w_h\|^2$.
Under the assumption on $\Delta t$, we have $\frac{1}{\Delta t} - \frac{b^2}{2\kappa_0} - s \ge \frac{1}{2 \Delta t}$, and hence
\begin{equation*}
    \frac{1}{2\Delta t} \|w_h\|^2 + \frac{\kappa_0}{2} \|\nabla w_h\|^2 \le \|f\| \|w_h\| + \|\mathbf{g}\| \|\nabla w_h\|.
\end{equation*}
Multiplying by $2 \Delta t$ and using $\|w_h\| \le \tnorm{w_h}$ and $\Delta t^{1/2} \|\nabla w_h\| \le \tnorm{w_h}$ gives \eqref{eq:energy}.
\end{proof}

Lemma~\ref{lem:energy} shows that a perturbation entering a step of the algorithm as a source term in $L^2(\Omega)$
is damped by a factor $\Delta t$, while a perturbation entering as a flux, i.e., through $\mathbf{g}$, is damped by $\Delta t^{1/2}$.
An important exception is a perturbation of an advection field, e.g., of $2 \nabla A_h / A_h$ in \eqref{eq:iterative-2}.
The corresponding flux contains $\nabla \delta A_h$, which the norm \eqref{eq:tnorm} controls only with weight $\Delta t^{1/2}$,
i.e., $\|\nabla \delta A_h\| \le \Delta t^{-1/2} \tnorm{\delta A_h}$.
Combined with \eqref{eq:energy}, the two factors cancel, so that the perturbation of the advection field
is transferred without damping: $\tnorm{\delta \rho_h} \le C \tnorm{\delta A_h}$.
This observation explains a difference between the fully decoupled and the coupled approaches, as we will see in Remark~\ref{rem:comparison}.

To apply Lemma~\ref{lem:energy} to the four steps in the evaluation of $G_d$, we need bounds on the coefficients of the corresponding problems and on their variations.
The advection fields $2 \nabla A_h / A_h$ and $2 \nabla H_h / H_h$ are bounded if $A_h$ and $H_h$ are bounded away from zero and have bounded gradients.
The reaction terms and the right-hand sides involve $\rho_h$, $\pi_h$, and $e^{-\pi_h}$, which are bounded if $\rho_h$ and $\pi_h$ are bounded and nonnegative.
In addition, the update of $H$ requires bounded gradients of $\rho_h$ and $\pi_h$ (see Lemma~\ref{lem:nodal} below).
Finally, these bounds have to hold along the whole iteration, i.e., in a set that is mapped into itself by $G_d$. This leads to the following assumption.

\begin{assump}\label{assump:bounds}
At time step $t^{n+1}$, there exist a closed convex set $\mathcal{K}_d \subset X$ with nonempty interior and positive constants
$a_{\min}$, $a_{\max}$, $r_{\max}$, $p_{\max}$, and $g_{\max}$ such that $G_d(\mathcal{K}_d) \subset \mathcal{K}_d$ and, for all $\mathbf{u}^{k} \in \mathcal{K}_d$, the functions computed in the evaluation of $G_d(\mathbf{u}^{k})$ satisfy
\begin{align*}
    & a_{\min} \le A^{k+1}_h \le a_{\max}, \quad a_{\min} \le H^{k+1}_h \le a_{\max}, \quad
    0 \le \rho^{k}_h \le r_{\max}, \quad 0 \le \pi^{k}_h \le p_{\max}, \\
    & |\nabla A^{k+1}_h|, \ |\nabla H^{k+1}_h|, \ |\nabla \rho^{k}_h|, \ |\nabla \pi^{k}_h| \le g_{\max} \quad \text{in } \overline{\Omega}.
\end{align*}
\end{assump}

The lower bounds on $\rho_h$ and $\pi_h$ in Assumption~\ref{assump:bounds} are discrete positivity properties.
For the Galerkin discretization of the $\rho$ and $\pi$ equations, this requires the local mesh P\'eclet numbers $h|\nabla A_h|/A_h$ and $h|\nabla H_h|/H_h$ to be sufficiently small and, because of the consistent mass and reaction matrices,
$\Delta t$ not to be too small compared with $h^2$.
Since the P\'eclet numbers depend on the discrete solution itself, these conditions cannot be guaranteed a priori, and the bounds remain an assumption on the computed solution.
In Sec.~\ref{sec:num_bounds}, we check the bounds in Assumption~\ref{assump:bounds} numerically.
In addition to Lemma~\ref{lem:energy}, we need a stability estimate for the nodal update of $H$ in \eqref{eq:iterative-3-v}.
We denote by $I_h$ the nodal interpolant onto $V_h(\Omega)$. Throughout the rest of this section, we assume that the mesh is shape regular and quasi-uniform.
The analysis below is carried out for $\mathbb{Q}_1$ elements, but it extends to higher-degree elements with constants depending on the degree.

\begin{lem}\label{lem:nodal}
Let $(\rho_h, A_h, \pi_h), (\tilde{\rho}_h, \tilde{A}_h, \tilde{\pi}_h) \in V_h(\Omega)^3$ satisfy the bounds
in Assumption~\ref{assump:bounds} on the corresponding variables and their gradients. Then, for $\Delta t \le 1$,
\begin{equation}\label{eq:nodal}
    \tnorm{\delta I_h (\rho_h A_h e^{-\pi_h})} \le C_{\rm mesh} \left( L + L_2 \, g_{\max} \right) \left( \tnorm{\delta \rho_h} + \tnorm{\delta A_h} + \tnorm{\delta \pi_h} \right),
\end{equation}
where $\delta I_h (\rho_h A_h e^{-\pi_h}) = I_h (\rho_h A_h e^{-\pi_h}) - I_h (\tilde{\rho}_h \tilde{A}_h e^{-\tilde{\pi}_h})$, $L = \max\{a_{\max}, r_{\max}, r_{\max} a_{\max}\}$, $L_2 = \max\{1, L\}$,
and $C_{\rm mesh}$ depends only on the shape-regularity and quasi-uniformity constants of the mesh.
\end{lem}

The proof is given in Appendix~\ref{app:nodal}.

\begin{prop}\label{prop:contraction}
Let Assumption~\ref{assump:bounds} hold and define
\begin{equation}\label{eq:dt0}
    \Delta t_0 = \min \left\{ 1, \ \frac{a_{\min}^2}{4 g_{\max}^2}, \ \frac{1}{2 \max\{0, r_{\max} - 1\}} \right\},
\end{equation}
where the last term is omitted if $r_{\max} \le 1$. Then, for all $\Delta t \le \Delta t_0$,
\begin{equation}\label{eq:contraction}
    \tnorm{G_d(\mathbf{u}) - G_d(\tilde{\mathbf{u}})} \le C_d \, \Delta t \, \tnorm{\mathbf{u} - \tilde{\mathbf{u}}} \quad \text{for all } \mathbf{u}, \tilde{\mathbf{u}} \in \mathcal{K}_d,
\end{equation}
where $C_d > 0$ depends only on the bounds in Assumption~\ref{assump:bounds}, on $\eta$, $\tau$, and $q$, and on the shape-regularity and quasi-uniformity constants of the mesh, but not on $\Delta t$.
Hence, for $\Delta t < \min\{\Delta t_0, C_d^{-1}\}$, $G_d$ is a contraction on $\mathcal{K}_d$ and has a unique fixed point $\mathbf{u}^* \in \mathcal{K}_d$,
i.e., the space- and time-discrete problem at time $t^{n+1}$ has a unique solution in $\mathcal{K}_d$, and
\begin{equation}\label{eq:k-step-d}
    \tnorm{\mathbf{u}^{k} - \mathbf{u}^*} \le (C_d \Delta t)^k \, \tnorm{\mathbf{u}^{0} - \mathbf{u}^*}, \quad k \ge 0.
\end{equation}
In particular, the Jacobian $J_d$ of $G_d$ at $\mathbf{u}^*$ satisfies $\tnorm{J_d^k} \le (C_d \Delta t)^k$, $k \ge 1$, and its spectral radius is bounded by $C_d \Delta t$.
\end{prop}

\begin{proof}
Let $\mathbf{u} = (\rho^{k}_h, \pi^{k}_h)$ and $\tilde{\mathbf{u}} = (\tilde{\rho}^{k}_h, \tilde{\pi}^{k}_h)$ be in $\mathcal{K}_d$.
Since $G_d(\mathcal{K}_d) \subset \mathcal{K}_d$, the bounds in Assumption~\ref{assump:bounds} hold for the inputs and for all the functions computed in the evaluation of $G_d$, including $\rho^{k+1}_h$ and $\pi^{k+1}_h$.
We will repeatedly use that, for $\pi, \tilde{\pi} \ge 0$, $e^{-\pi} \le 1$ and $|e^{-\pi} - e^{-\tilde{\pi}}| \le |\pi - \tilde{\pi}|$, and that
$\tnorm{\delta \rho^{k}_h} + \tnorm{\delta \pi^{k}_h} \le \sqrt{2} \, \tnorm{\delta \mathbf{u}}$.
Below, $K_A$, $K_{\rho A}$, $K_\rho$, $K_H$, and $K_\pi$ denote positive constants that depend only on the bounds in Assumption~\ref{assump:bounds}, on $\eta$ and $q$,
and on the constants $C_{\rm mesh}$, $L$, and $L_2$ in Lemma~\ref{lem:nodal}.

\emph{Step 1.} The difference $\delta A_h = A^{k+1}_h - \tilde{A}^{k+1}_h$ solves \eqref{eq:model-problem} with $\kappa_0 = \eta$, $\boldsymbol{\beta} = \mathbf{0}$, $\sigma = 1 - c^{k}_h$, $\mathbf{g} = \mathbf{0}$,
and $f = \delta c_h \, \tilde{A}^{k+1}_h$, where we recall $c^{k}_h = \rho^{k}_h e^{-\pi^{k}_h}$. %
Here, $s = \max\{0, r_{\max} - 1\}$ since $0 \le c^{k}_h \le r_{\max}$. Since $b = 0$, the condition of Lemma~\ref{lem:energy} reads
$\Delta t \le 1 / (2 \max\{0, r_{\max} - 1\})$, which holds for $\Delta t \le \Delta t_0$.
This is the origin of the last term in \eqref{eq:dt0}.
This condition is closely related to the sufficient condition for the matrix in \eqref{eq:mat1bis} to be symmetric positive definite.
Since $\delta c_h = \delta \rho^{k}_h \, e^{-\pi^{k}_h} + \tilde{\rho}^{k}_h (e^{-\pi^{k}_h} - e^{-\tilde{\pi}^{k}_h})$, we have $\|\delta c_h\| \le \|\delta \rho^{k}_h\| + r_{\max} \|\delta \pi^{k}_h\|$, and \eqref{eq:energy} gives
\begin{equation}\label{eq:step1}
    \tnorm{\delta A_h} \le \frac{2 a_{\max}}{\min\{1, \eta\}} \Delta t \left( \|\delta \rho^{k}_h\| + r_{\max} \|\delta \pi^{k}_h\| \right) \le K_A \Delta t \, \tnorm{\delta \mathbf{u}}.
\end{equation}

\emph{Step 2.} The difference $\delta \rho_h = \rho^{k+1}_h - \tilde{\rho}^{k+1}_h$ solves \eqref{eq:model-problem} with $\kappa_0 = 1$, $\boldsymbol{\beta} = 2 \nabla A^{k+1}_h / A^{k+1}_h$,
so that $b = 2 g_{\max} / a_{\min}$, $\sigma = A^{k+1}_h e^{-\pi^{k}_h} \ge 0$, so that $s = 0$, and
\begin{equation*}
    \mathbf{g} = 2 \tilde{\rho}^{k+1}_h \, \delta \left( \frac{\nabla A^{k+1}_h}{A^{k+1}_h} \right), \qquad f = - \delta (A^{k+1}_h e^{-\pi^{k}_h}) \, \tilde{\rho}^{k+1}_h + q \, \delta (e^{-\pi^{k}_h}).
\end{equation*}
The condition of Lemma~\ref{lem:energy} reads $\Delta t \le a_{\min}^2 / (4 g_{\max}^2)$, which holds for $\Delta t \le \Delta t_0$. Note that, dropping the superscript $k+1$ for brevity,
\begin{equation*}
    \delta \left( \frac{\nabla A_h}{A_h} \right) = \frac{\nabla A_h}{A_h} - \frac{\nabla \tilde{A}_h}{\tilde{A}_h} = \frac{\nabla \delta A_h}{A_h} - \frac{\nabla \tilde{A}_h \, \delta A_h}{A_h \tilde{A}_h}.
\end{equation*}
Given this,
the bounds in Assumption~\ref{assump:bounds}, and $\Delta t^{1/2} \|\nabla \delta A_h\| \le \tnorm{\delta A_h}$, for $\Delta t \le 1$ we have
\begin{equation*}
    2 \Delta t^{1/2} \|\mathbf{g}\| \le 4 r_{\max} \left( \frac{1}{a_{\min}} + \frac{g_{\max}}{a_{\min}^{2}} \right) \tnorm{\delta A_h}.
\end{equation*}
Moreover, $|\delta (A^{k+1}_h e^{-\pi^{k}_h})| \le |\delta A_h| + a_{\max} |\delta \pi^{k}_h|$, so that
\begin{equation*}
    2 \Delta t \|f\| \le 2 r_{\max} \tnorm{\delta A_h} + 2 \Delta t \, (r_{\max} a_{\max} + q) \tnorm{\delta \pi^{k}_h}.
\end{equation*}
Then, \eqref{eq:energy} and \eqref{eq:step1} give
\begin{equation}\label{eq:step2}
    \tnorm{\delta \rho_h} \le K_{\rho A} \tnorm{\delta A_h} + 2 \Delta t \, (r_{\max} a_{\max} + q) \tnorm{\delta \pi^{k}_h} \le K_\rho \Delta t \, \tnorm{\delta \mathbf{u}}.
\end{equation}

\emph{Step 3.} From \eqref{eq:iterative-3-v}, $\delta H_h = \frac{\Delta t}{\tau + \Delta t} \, \delta I_h (\rho^{k+1}_h A^{k+1}_h e^{-\pi^{k}_h})$.
By Assumption~\ref{assump:bounds}, the functions $\rho^{k+1}_h$, $A^{k+1}_h$, and $\pi^{k}_h$ and their counterparts computed from $\tilde{\mathbf{u}}$ satisfy the hypotheses of Lemma~\ref{lem:nodal}. Hence, by \eqref{eq:nodal}, \eqref{eq:step1}, \eqref{eq:step2}, and $\Delta t \le 1$,
\begin{equation}\label{eq:step3}
    \tnorm{\delta H_h} \le \frac{\Delta t}{\tau} C_{\rm mesh} (L + L_2 \, g_{\max}) \left( \tnorm{\delta \rho_h} + \tnorm{\delta A_h} + \tnorm{\delta \pi^{k}_h} \right) \le K_H \frac{\Delta t}{\tau} \tnorm{\delta \mathbf{u}}.
\end{equation}

\emph{Step 4.} The difference $\delta \pi_h = \pi^{k+1}_h - \tilde{\pi}^{k+1}_h$ solves \eqref{eq:model-problem} with $\kappa_0 = 1$, $\boldsymbol{\beta} = 2 \nabla H^{k+1}_h / H^{k+1}_h$, so that $b = 2 g_{\max} / a_{\min}$,
$\sigma = 0$, $f = 0$, and $\mathbf{g} = 2 \tilde{\pi}^{k+1}_h \, \delta (\nabla H^{k+1}_h / H^{k+1}_h)$. Proceeding as in Step 2 with $\tilde{\pi}^{k+1}_h \le p_{\max}$, we obtain
\begin{equation}\label{eq:step4}
    \tnorm{\delta \pi_h} %
    \le K_\pi \frac{\Delta t}{\tau} \tnorm{\delta \mathbf{u}}.
\end{equation}

Finally, $\tnorm{G_d(\mathbf{u}) - G_d(\tilde{\mathbf{u}})} \le \tnorm{\delta \rho_h} + \tnorm{\delta \pi_h}$, and \eqref{eq:step2} and \eqref{eq:step4} give \eqref{eq:contraction} with $C_d = K_\rho + K_\pi / \tau$.
Estimate \eqref{eq:k-step-d} follows from the Banach fixed-point theorem.

Finally, $G_d$ is continuously differentiable on $\mathcal{K}_d$ (see Remark~\ref{rem:AA-norm}). Hence, \eqref{eq:contraction} implies $\tnorm{J_d} \le C_d \Delta t$ at the interior points of $\mathcal{K}_d$ and, since $\mathcal{K}_d$ is convex with nonempty interior, by continuity on all of $\mathcal{K}_d$.
This gives the bounds on $J_d^k$ and on its spectral radius.
\end{proof}

Proposition~\ref{prop:contraction} implies that the number of iterations per time step is expected to grow with $\Delta t$,
and that the fully decoupled iteration may fail to converge when $C_d \Delta t \geq 1$.

\begin{remark}\label{rem:AA-norm}
The map $G_d$ is a composition of smooth maps: products, exponentials, $1/A_h$ and $1/H_h$ with $A_h, H_h \ge a_{\min}$, and inverses of the matrices in \eqref{eq:mat1bis}, \eqref{eq:mat2bis}, and \eqref{eq:mat4},
which are uniformly nonsingular under the assumptions of Proposition~\ref{prop:contraction}. Hence, $G_d$ is continuously differentiable on $\mathcal{K}_d$ with Lipschitz continuous derivative.
Together with \eqref{eq:contraction}, this implies that AA($m$) converges locally R-linearly \cite{TothKelley2015} and that it improves the linear convergence rate $C_d \Delta t$ by the gain factor of the
least-squares problem, up to higher-order terms \cite{Evans2020}. The latter result requires the least-squares problem in \eqref{eq:AA-update} to be posed in the inner-product norm in which $G_d$ is contractive.
With $B$ given by \eqref{eq:B}, we have $\|\mathbf{x}\|_B^2 = \sum_j w_j^2 \tnorm{x_j}^2$, where $x_j$ is the $j$-th block of $\mathbf{x}$, i.e., $\|\cdot\|_B$ is equivalent to \eqref{eq:tnorm} with block weights.
In this norm, $G_d$ is a contraction with constant $(\max_j w_j / \min_j w_j) \, C_d \Delta t$, which is still $\mathcal{O}(\Delta t)$. Hence, the requirement is met.
When this constant is larger than one, these results no longer apply, although AA may still converge. See \cite{pollock2021anderson}.
\end{remark}

\section{A coupled approach}\label{sec:coupled}

In the fully decoupled approach of Sec.~\ref{sec:decoupled}, the coupling between $A$ and $\rho$ through the reaction term $\rho A e^{-\pi}$
is treated in a staggered way: in Step 1, $\rho$ is taken from the previous iteration,
while in Step 2, $A$ is taken from the current one. When this coupling is strong, the fixed-point iteration can converge slowly or not at all.
Here, we consider an alternative approach that solves for $A$ and $\rho$ simultaneously. The product $\rho A$ is linearized
with Newton's method, i.e., by its first-order Taylor expansion around $(\rho^{k}_h, A^{k}_h)$,
\begin{equation*}
    \rho^{k+1}_h A^{k+1}_h \approx \rho^{k+1}_h A^{k}_h + \rho^{k}_h A^{k+1}_h - \rho^{k}_h A^{k}_h,
\end{equation*}
consistently in both equations, while the advection term in the equation for $\rho$ is evaluated with $A^{k}_h$
and $e^{-\pi}$ is evaluated with $\pi^{k}_h$, so that the resulting problem is linear. Variables $H$ and $\pi$ are then updated as in Sec.~\ref{sec:decoupled}.
As in Sec.~\ref{sec:decoupled}, we directly state the algorithm in its fully discretized form.

The algorithm is as follows:
at every time step $t^{n+1}$, given $A^n_h$, $\rho^n_h$, $\pi^n_h$, and $H^n_h$, set
$A^{0}_h = A^n_h$, $\rho^{0}_h = \rho^{n}_h$, and $\pi^{0}_h = \pi^{n}_h$ and
perform at iteration $k+1$, $k \ge 0$ (see notation explanation in Remark~\ref{rem:notation}):
\begin{itemize}
    \item[-] \emph{Step 1}: find $(A^{k+1}_h, \rho^{k+1}_h) \in V_h(\Omega)^2$ such that
    \begin{align}
           &\frac{1}{\Delta t}\int_\Omega A^{k+1}_h u_h\;d\Omega +
           \int_\Omega \eta\nabla A^{k+1}_h \cdot \nabla u_h\;d\Omega +\int_\Omega A^{k+1}_h u_h \;d\Omega \cl
           & \quad \quad -\int_\Omega e^{- \pi^{k}_h} \left( \rho^{k+1}_h A^{k}_h + \rho^{k}_h A^{k+1}_h \right) u_h \;d\Omega = \int_\Omega \Ast  u_h \;d\Omega  \cl
           & \quad \quad + \frac{1}{\Delta t}\int_\Omega A^{n}_h u_h\;d\Omega
           - \int_\Omega \rho^{k}_h A^{k}_h e^{- \pi^{k}_h} u_h \;d\Omega, \label{eq:coupled-1} \\
           & \frac{1}{\Delta t} \int_\Omega \rho^{k+1}_h v_h \;d\Omega +\int_\Omega  \Big(\nabla\rho^{k+1}_h - \frac{2\nabla A^{k}_h}{A^{k}_h} \rho^{k+1}_h\Big) \cdot \nabla v_h\;d\Omega \cl
           & \quad \quad + \int_\Omega e^{- \pi^{k}_h} \left( \rho^{k+1}_h A^{k}_h + \rho^{k}_h A^{k+1}_h \right) v_h \;d\Omega = \int_\Omega q e^{-\pi^{k}_h} v_h\;d\Omega  \cl
           & \quad \quad + \frac{1}{\Delta t} \int_\Omega \rho^{n}_h v_h \;d\Omega
           + \int_\Omega \rho^{k}_h A^{k}_h e^{- \pi^{k}_h} v_h \;d\Omega, \label{eq:coupled-2}
    \end{align}
    for all $(u_h, v_h) \in V_h(\Omega)^2$.
    \item[-] \emph{Step 2}: set $H^{k+1}_h \in V_h(\Omega)$ such that \eqref{eq:iterative-3} holds at the nodes of $\mathcal{T}_h$.
    \item[-] \emph{Step 3}: find $\pi^{k+1}_h \in V_h(\Omega)$ such that \eqref{eq:iterative-4} holds for all $w_h \in V_h(\Omega)$.
    \item[-] \emph{Step 4}: Check stopping criterion \eqref{eq:crit}. If not satisfied, repeat steps 1-4.
    If satisfied, set $A^{n+1}_h = A^{k+1}_h$, $\rho^{n+1}_h = \rho^{k+1}_h$, $\pi^{n+1}_h = \pi^{k+1}_h$, and $H^{n+1}_h = H^{k+1}_h$.
\end{itemize}
Since the linearization in \eqref{eq:coupled-1}-\eqref{eq:coupled-2} is consistent,
the fixed points of the coupled approach, as those of the fully decoupled one,
are solutions of the spatial discretization of the time-discrete problem \eqref{eq:implicit-1}-\eqref{eq:implicit-4}.
Hence, both approaches approximate the same discrete solution, and their computational costs can be compared directly.

Let $c^{k}_h = \rho^{k}_h e^{-\pi^{k}_h}$, as in \eqref{eq:mat1bis}, and $a^{k}_h = A^{k}_h e^{-\pi^{k}_h}$.
Then, the matrix form of Step 1 is given by: find the vectors $\bA^{k+1}, \br^{k+1} \in \mathbb{R}^{N_Q}$ such that
\begin{equation}\label{eq:mat-coupled}
    \begin{bmatrix}
        \dfrac{1 + \Delta t}{\Delta t} M + \eta K - N(c^{k}_h) & - N(a^{k}_h) \\[2ex]
        N(c^{k}_h) & \dfrac{1}{\Delta t} M + K - D(A^{k}_h) + N(a^{k}_h)
    \end{bmatrix}
    \begin{bmatrix}
        \bA^{k+1} \\[2ex] \br^{k+1}
    \end{bmatrix}
    =
    \begin{bmatrix}
        \bb_A^{k+1} \\[2ex] \bb_\rho^{k+1}
    \end{bmatrix},
\end{equation}
where
\begin{equation*}
    \bb_A^{k+1} = M  \bA^{\mathrm{st}} + \frac{1}{\Delta t} M \bA^{n} - N(c^{k}_h) \bA^{k},
\end{equation*}
and $\bb_\rho^{k+1}$ has $i$-th component
\begin{equation*}
    \bb_{\rho,i}^{k+1} = \int_\Omega q e^{-\pi^{k}_h} \phi_i \;d\Omega + \frac{1}{\Delta t} (M\br^{n})_i + \left( N(c^{k}_h) \bA^{k} \right)_i,
\end{equation*}
for $i = 1, \dots, N_Q$. Steps 2 and 3 are performed as in Sec.~\ref{sec:decoupled}, i.e., with \eqref{eq:iterative-3-v} and \eqref{eq:mat4},
where now $\bd^{k+1}$ is computed from $\bA^{k+1}$ and $\br^{k+1}$ given by \eqref{eq:mat-coupled}.

Compared to the fully decoupled approach, the coupled approach requires the solution of one linear system of size $2N_Q$,
i.e., \eqref{eq:mat-coupled}, and one of size $N_Q$, i.e., \eqref{eq:mat4}, per iteration, instead of three linear systems of size $N_Q$.
The $(1,1)$ block of the matrix in \eqref{eq:mat-coupled}  coincides with the matrix in \eqref{eq:mat1bis}, while the $(2,2)$ block differs
from the matrix in \eqref{eq:mat2bis} only in that $A^{k}_h$ replaces $A^{k+1}_h$ in $D(\cdot)$ and $N(\cdot)$.
The off-diagonal blocks are weighted mass matrices that account for the Newton linearization of the reaction term.
Since the off-diagonal blocks differ and $D(A^{k}_h)$ is nonsymmetric, the matrix in \eqref{eq:mat-coupled} is nonsymmetric.

\subsection{Fixed-point formulation and Anderson acceleration}
\label{sec:AA-coupled}

Given the solution at time $t^n$, iteration $k+1$ of the coupled approach depends on $\bA^{k}$, $\br^{k}$, and $\bp^{k}$:
$\bA^{k}$ and $\br^{k}$ enter the linearized reaction terms in \eqref{eq:mat-coupled}, $\bA^{k}$ also enters the advection term through $D(A^{k}_h)$,
and $\bp^{k}$ enters through $e^{-\pi^{k}_h}$. Instead, $\bH^{k}$ plays no role.
Let $\mathbf{u} = (\bA, \br, \bp) \in \mathbb{R}^{3N_Q}$ and let $G_c: \mathbb{R}^{3N_Q} \to \mathbb{R}^{3N_Q}$ be the map
that takes $\mathbf{u}^{k}$ to $\mathbf{u}^{k+1}$ through one iteration of the coupled approach.
Then, the coupled algorithm takes the form \eqref{eq:fixed-point} with $G = G_c$ and $\mathbf{u}^{0} = (\bA^{n}, \br^{n}, \bp^{n})$.
Each evaluation of $G_c$ consists of the following operations:
\begin{enumerate}
    \item assemble $N(c^{k}_h)$, $N(a^{k}_h)$, and $D(A^{k}_h)$, with $c^{k}_h = \rho^{k}_h e^{-\pi^{k}_h}$ and $a^{k}_h = A^{k}_h e^{-\pi^{k}_h}$, and solve \eqref{eq:mat-coupled} for $(\widetilde{\bA}^{k+1}, \widetilde{\br}^{k+1})$;
    \item compute $\bd^{k+1}$ from $\widetilde{\bA}^{k+1}$, $\widetilde{\br}^{k+1}$, and $\bp^{k}$, and compute $\bH^{k+1}$ with \eqref{eq:iterative-3-v};
    \item solve \eqref{eq:mat4} for $\widetilde{\bp}^{k+1}$, and set $G_c(\mathbf{u}^{k}) = (\widetilde{\bA}^{k+1}, \widetilde{\br}^{k+1}, \widetilde{\bp}^{k+1})$.
\end{enumerate}
The coupled algorithm with AA($m$) is given by Alg.~\ref{alg:AA-generic} with $G = G_c$, where the stopping criterion is replaced by \eqref{eq:crit} applied to $G_c(\mathbf{u}^{k})$ and $\mathbf{u}^{k}$.
Compared to the fully decoupled approach, the iterate has size $3N_Q$ instead of $2N_Q$, which slightly increases the cost of the least-squares problem in \eqref{eq:AA-update}.

The considerations in Remark~\ref{rem:AA-positivity} apply also to the coupled approach, with one important difference.
In the fully decoupled approach, the attractiveness is always obtained from a linear solve, i.e., \eqref{eq:mat1bis}.
In the coupled approach, the accelerated attractiveness $\bA^{k+1}$ is part of the new iterate and enters the next evaluation of $G_c$ directly,
through the advection field $2 \nabla A^{k+1}_h / A^{k+1}_h$ in $D(A^{k+1}_h)$. Hence, the positivity of the accelerated attractiveness is essential,
and the safeguard described in Remark~\ref{rem:AA-positivity} should include $\bA^{k+1}$.

\subsection{Convergence of the coupled iteration}\label{sec:contraction-coupled}

In this section, we derive the counterpart of Proposition~\ref{prop:contraction} for the map $G_c$. We identify $\mathbf{u} = (\bA, \br, \bp)$ with $(A_h, \rho_h, \pi_h) \in V_h(\Omega)^3$
and use the norm $\tnorm{(x_A, x_\rho, x_\pi)}_1 = \tnorm{x_A} + \tnorm{x_\rho} + \tnorm{x_\pi}$, with $\tnorm{\cdot}$ defined in \eqref{eq:tnorm}.
For a linear map $J: V_h(\Omega)^3 \to V_h(\Omega)^3$, e.g., the Jacobian of $G_c$, we also write
$\tnorm{J}_1 = \sup_{\mathbf{x} \neq \mathbf{0}} \tnorm{J \mathbf{x}}_1 / \tnorm{\mathbf{x}}_1$.
To compare the two approaches in the same space, we note that $G_d$ can be viewed as a map on $V_h(\Omega)^3$ that does not depend on its first argument,
so that Proposition~\ref{prop:contraction} also holds in this norm, since its proof bounds $\tnorm{\delta \rho_h} + \tnorm{\delta \pi_h}$ and $\tnorm{\cdot} \le \tnorm{\cdot}_1$.
As in Sec.~\ref{sec:contraction}, we assume that the mesh is shape regular and quasi-uniform.
Since the input attractiveness $A^{k}_h$ now enters the advection term in \eqref{eq:coupled-2}, it must satisfy the same bounds as $A^{k+1}_h$.

\begin{assump}\label{assump:bounds-coupled}
At time step $t^{n+1}$, there exists a closed convex set $\mathcal{K}_c \subset V_h(\Omega)^3$ with nonempty interior such that $G_c(\mathcal{K}_c) \subset \mathcal{K}_c$
and, for all $\mathbf{u}^{k} \in \mathcal{K}_c$, the bounds in Assumption~\ref{assump:bounds} hold for the functions
computed in the evaluation of $G_c(\mathbf{u}^{k})$ and, in addition, $a_{\min} \le A^{k}_h \le a_{\max}$ and
$|\nabla A^{k}_h| \le g_{\max}$ in $\overline{\Omega}$.
\end{assump}

As for Assumption~\ref{assump:bounds}, the lower bound on $\rho_h$ in Assumption~\ref{assump:bounds-coupled} is a discrete positivity property, with the mesh P\'eclet number for the $\rho$ equation now computed from the lagged field $A^{k}_h$.
However, the Newton linearization of the reaction term introduces in the $\rho$ equation a coupling with $A^{k+1}_h$ with positive off-diagonal entries,
so the matrix of the coupled system in \eqref{eq:mat-coupled} is not an M-matrix, and positivity of $\rho_h$ cannot be guaranteed even for small P\'eclet numbers.
Nevertheless, $\rho_h$ remains positive in all our computations, as we will see in Sec.~\ref{sec:num_bounds}.

\begin{prop}\label{prop:contraction-coupled}
Let Assumption~\ref{assump:bounds-coupled} hold and define
\begin{equation*}
    \Delta t_0 = \min \left\{ 1, \ \left( \frac{4 g_{\max}^2}{a_{\min}^2} + 3 r_{\max} + a_{\max} \right)^{-1} \right\}.
\end{equation*}
Then, there exist constants $C_c > 0$ and $C_{\rm adv} > 0$, which depend only on the bounds in Assumption~\ref{assump:bounds-coupled}, on $\eta$, $\tau$, and $q$,
and on the shape-regularity and quasi-uniformity constants of the mesh, but not on $\Delta t$, such that for all $\Delta t \le \Delta t_0$ and all $\mathbf{u}, \tilde{\mathbf{u}} \in \mathcal{K}_c$
\begin{align}
    \tnorm{\delta A^{k+1}_h} + \tnorm{\delta \pi^{k+1}_h} & \le C_c \Delta t \, \tnorm{\mathbf{u} - \tilde{\mathbf{u}}}_1, \label{eq:lip-c-1} \\
    \tnorm{\delta \rho^{k+1}_h} & \le C_{\rm adv} \tnorm{\delta A^{k}_h} + C_c \Delta t \, \tnorm{\mathbf{u} - \tilde{\mathbf{u}}}_1. \label{eq:lip-c-2}
\end{align}
Moreover, if $\mathbf{u}^* \in \mathcal{K}_c$ is a fixed point of $G_c$, the Jacobian $J_c$ of $G_c$ at $\mathbf{u}^*$ has vanishing columns associated with $\rho$ and,
for $C_c \Delta t \le 1/2$,
\begin{equation}\label{eq:k-step-c}
    \tnorm{J_c^k}_1 \le (C_{\rm adv} + 1) (C_c \Delta t)^{k-1}, \quad k \ge 1.
\end{equation}
In particular, the spectral radius of $J_c$ is bounded by $C_c \Delta t$.
\end{prop}

\begin{proof}
Let $\mathbf{u} = (A^{k}_h, \rho^{k}_h, \pi^{k}_h)$ and $\tilde{\mathbf{u}} = (\tilde{A}^{k}_h, \tilde{\rho}^{k}_h, \tilde{\pi}^{k}_h)$ be in $\mathcal{K}_c$.
Since $G_c(\mathcal{K}_c) \subset \mathcal{K}_c$, the bounds in Assumption~\ref{assump:bounds-coupled} hold for the inputs and for all the functions computed in the evaluation of $G_c$.
Below, $K_S$, $K_\rho$, $K_A$, $K_H$, and $K_\pi$ denote positive constants that depend only on the bounds in Assumption~\ref{assump:bounds-coupled}, on $\eta$ and $q$,
and on the constants $C_{\rm mesh}$, $L$, and $L_2$ in Lemma~\ref{lem:nodal}. We use the same elementary inequalities as in the proof of Proposition~\ref{prop:contraction}.

Let $\boldsymbol{\alpha}^{k} = 2 \nabla A^{k}_h / A^{k}_h$ denote the advection field in \eqref{eq:coupled-2}.
Subtracting \eqref{eq:coupled-1}-\eqref{eq:coupled-2} written for $\mathbf{u}$ and $\tilde{\mathbf{u}}$,
we find that $(\delta A_h, \delta \rho_h) = (\delta A^{k+1}_h, \delta \rho^{k+1}_h)$ satisfies, for all $u_h, v_h \in V_h(\Omega)$,
\begin{align}
    & \frac{1}{\Delta t} (\delta A_h, u_h) + \eta (\nabla \delta A_h, \nabla u_h) + \big( (1 - e^{-\pi^{k}_h} \tilde{\rho}^{k}_h) \delta A_h, u_h \big) \cl
    & \quad \quad - (e^{-\pi^{k}_h} A^{k}_h \delta \rho_h, u_h) = (S, u_h), \label{eq:diff-c-1} \\
    & \frac{1}{\Delta t} (\delta \rho_h, v_h) + (\nabla \delta \rho_h, \nabla v_h) - (\boldsymbol{\alpha}^{k} \delta \rho_h, \nabla v_h) + (e^{-\pi^{k}_h} A^{k}_h \delta \rho_h, v_h) + (e^{-\pi^{k}_h} \tilde{\rho}^{k}_h \delta A_h, v_h) \cl
    & \quad \quad = - (S, v_h) + (q \, \delta (e^{-\pi^{k}_h}), v_h) + (\tilde{\rho}^{k+1}_h \, \delta \boldsymbol{\alpha}^{k}, \nabla v_h), \label{eq:diff-c-2}
\end{align}
where
\begin{align*}
    S = & \ e^{-\pi^{k}_h} \left[ (\tilde{\rho}^{k+1}_h - \tilde{\rho}^{k}_h) \, \delta A^{k}_h + (A^{k+1}_h - A^{k}_h) \, \delta \rho^{k}_h \right] \\
    & + \delta (e^{-\pi^{k}_h}) \left[ \tilde{\rho}^{k+1}_h \tilde{A}^{k}_h + \tilde{\rho}^{k}_h (\tilde{A}^{k+1}_h - \tilde{A}^{k}_h) \right]
\end{align*}
collects the variations of the linearized reaction terms with respect to the input. Since $|\tilde{\rho}^{k+1}_h - \tilde{\rho}^{k}_h| \le r_{\max}$, $|A^{k+1}_h - A^{k}_h| \le a_{\max}$, we have
\begin{equation}\label{eq:S-bound}
    \|S\| \le r_{\max} \|\delta A^{k}_h\| + a_{\max} \|\delta \rho^{k}_h\| + 2 r_{\max} a_{\max} \|\delta \pi^{k}_h\| \le K_S  \, \tnorm{\delta \mathbf{u}}_1.
\end{equation}
Note that the first two terms in $S$ are proportional to the update $\mathbf{u}^{k+1} - \mathbf{u}^{k}$ and vanish at a fixed point.

We first estimate $\delta \rho_h$. We take $u_h = \delta A_h$ in \eqref{eq:diff-c-1} and $v_h = \delta \rho_h$ in \eqref{eq:diff-c-2}, and sum the two resulting equations:
\begin{align}
    & \frac{1}{\Delta t} \left( \|\delta A_h\|^2 + \|\delta \rho_h\|^2 \right) + \eta \|\nabla \delta A_h\|^2 + \|\nabla \delta \rho_h\|^2
      + \big( (1 - e^{-\pi^{k}_h} \tilde{\rho}^{k}_h) \delta A_h, \delta A_h \big) \cl
    & \quad + (e^{-\pi^{k}_h} A^{k}_h \delta \rho_h, \delta \rho_h) - (e^{-\pi^{k}_h} A^{k}_h \delta \rho_h, \delta A_h) + (e^{-\pi^{k}_h} \tilde{\rho}^{k}_h \delta A_h, \delta \rho_h) - (\boldsymbol{\alpha}^{k} \delta \rho_h, \nabla \delta \rho_h) \cl
    & \quad = (S, \delta A_h - \delta \rho_h) + (q \, \delta (e^{-\pi^{k}_h}), \delta \rho_h) + (\tilde{\rho}^{k+1}_h \, \delta \boldsymbol{\alpha}^{k}, \nabla \delta \rho_h). \label{eq:sum}
\end{align}
We have
\begin{align*}
    \big( (1 - e^{-\pi^{k}_h} \tilde{\rho}^{k}_h) \delta A_h, \delta A_h \big) & \ge (1 - r_{\max}) \|\delta A_h\|^2, \\
    |(e^{-\pi^{k}_h} A^{k}_h \delta \rho_h, \delta A_h)| + |(e^{-\pi^{k}_h} \tilde{\rho}^{k}_h \delta A_h, \delta \rho_h)| & \le \frac{a_{\max} + r_{\max}}{2} \left( \|\delta A_h\|^2 + \|\delta \rho_h\|^2 \right), \\
    |(\boldsymbol{\alpha}^{k} \delta \rho_h, \nabla \delta \rho_h)| & \le \frac{1}{2} \|\nabla \delta \rho_h\|^2 + \frac{2 g_{\max}^2}{a_{\min}^2} \|\delta \rho_h\|^2, \\
    (e^{-\pi^{k}_h} A^{k}_h \delta \rho_h, \delta \rho_h) & \ge 0.
\end{align*}
For $\Delta t \le \Delta t_0$, the left-hand side of \eqref{eq:sum} is then bounded from below by
\begin{equation}\label{eq:bound}
    \frac{1}{2 \Delta t} \left( \|\delta A_h\|^2 + \|\delta \rho_h\|^2 \right) + \min\{\eta, \tfrac{1}{2}\} \left( \|\nabla \delta A_h\|^2 + \|\nabla \delta \rho_h\|^2 \right).
\end{equation}
Indeed, the bounds above reduce the coefficients of $\|\delta A_h\|^2$ and $\|\delta \rho_h\|^2$ to at least
$\frac{1}{\Delta t} - r_{\max} - \frac{a_{\max} + r_{\max}}{2}$ and $\frac{1}{\Delta t} - \frac{a_{\max} + r_{\max}}{2} - \frac{2 g_{\max}^2}{a_{\min}^2}$,
respectively, which are both larger than $\frac{1}{2 \Delta t}$ for $\Delta t \le \Delta t_0$, while the coefficient
of $\|\nabla \delta \rho_h\|^2$ is reduced to $\frac{1}{2}$.
Hence, the expression in \eqref{eq:bound} is bounded by the right-hand side of \eqref{eq:sum}. Multiplying \eqref{eq:bound} by $2 \Delta t$
and bounding the right-hand side of \eqref{eq:sum} with the Cauchy--Schwarz inequality, \eqref{eq:S-bound}, and  $\Delta t^{1/2} \|\delta \boldsymbol{\alpha}^{k}\| \le 2 \left( a_{\min}^{-1} + g_{\max} a_{\min}^{-2} \right) \tnorm{\delta A^{k}_h}$ (see proof of Proposition~\ref{prop:contraction}), we obtain
\begin{equation}\label{eq:step-c-rho}
    \tnorm{\delta \rho_h} \le C_{\rm adv} \tnorm{\delta A^{k}_h} + K_\rho \Delta t  \, \tnorm{\delta \mathbf{u}}_1, \qquad
    C_{\rm adv} = \frac{4 r_{\max}}{\min\{1, 2 \eta\}} \left( \frac{1}{a_{\min}} + \frac{g_{\max}}{a_{\min}^2} \right).
\end{equation}
Note that the first term, which comes from the advection field and contains the input $\delta A^{k}_h$, is not multiplied by a positive power of $\Delta t$.

Next, we estimate $\delta A_h$. Equation \eqref{eq:diff-c-1} is of the form \eqref{eq:model-problem} with $\kappa_0 = \eta$, $\boldsymbol{\beta} = \mathbf{0}$, $\sigma = 1 - e^{-\pi^{k}_h} \tilde{\rho}^{k}_h \ge 1 - r_{\max}$,
$\mathbf{g} = \mathbf{0}$, and $f = S + e^{-\pi^{k}_h} A^{k}_h \delta \rho_h$. The condition of Lemma~\ref{lem:energy} holds for $\Delta t \le \Delta t_0$, and \eqref{eq:energy}, \eqref{eq:S-bound}, and \eqref{eq:step-c-rho} give
\begin{equation}\label{eq:step-c-A}
    \tnorm{\delta A_h} \le \frac{2}{\min\{1, \eta\}} \Delta t \left( \|S\| + a_{\max} \|\delta \rho_h\| \right) \le K_A \Delta t  \, \tnorm{\delta \mathbf{u}}_1.
\end{equation}
Hence, the undamped contribution of $\delta A^{k}_h$ to $\delta \rho_h$ enters $\delta A_h$ only through the reaction term, i.e., as an $L^2$ source, and is damped by $\Delta t$.

The updates of $H$ and $\pi$ are the same as in the fully decoupled approach.
Hence, Steps 3 and 4 in the proof of Proposition~\ref{prop:contraction}, together with \eqref{eq:step-c-rho} and \eqref{eq:step-c-A}, give
\begin{equation}\label{eq:step-c-pi}
    \tnorm{\delta H^{k+1}_h} \le K_H \frac{\Delta t}{\tau} \tnorm{\delta \mathbf{u}}_1, \qquad \tnorm{\delta \pi^{k+1}_h} \le K_\pi \frac{\Delta t}{\tau} \tnorm{\delta \mathbf{u}}_1.
\end{equation}
Estimates \eqref{eq:step-c-rho}, \eqref{eq:step-c-A}, and \eqref{eq:step-c-pi} give \eqref{eq:lip-c-1}-\eqref{eq:lip-c-2} with $C_c = \max\{K_A + K_\pi / \tau, \, K_\rho\}$.

Finally, we consider the Jacobian of $G_c$. The map $G_c$ is continuously differentiable on $\mathcal{K}_c$, for the same reasons as $G_d$ (see Remark~\ref{rem:AA-norm}).
The input $\rho^{k}_h$ appears in the evaluation of $G_c$ only in the linearized reaction terms, i.e., in $e^{-\pi^{k}_h} \rho^{k}_h (A^{k+1}_h - A^{k}_h)$,
while the updates of $H$ and $\pi$ use the output $\rho^{k+1}_h$. Accordingly, $\delta \rho^{k}_h$ enters \eqref{eq:diff-c-1}-\eqref{eq:diff-c-2} only through
the term $e^{-\pi^{k}_h} (A^{k+1}_h - A^{k}_h) \, \delta \rho^{k}_h$ in $S$. At a fixed point, $A^{k+1}_h = A^{k}_h$, so this term vanishes:
$G_c$ does not depend on $\rho^{k}_h$ to first order, i.e., the columns of $J_c$ associated with $\rho$ vanish. Ordering the unknowns as $(A, \rho, \pi)$, $J_c$ has the block structure
\begin{equation*}
    J_c = \begin{bmatrix} J_{AA} & 0 & J_{A\pi} \\ J_{\rho A} & 0 & J_{\rho \pi} \\ J_{\pi A} & 0 & J_{\pi \pi} \end{bmatrix},
\end{equation*}
and, by \eqref{eq:lip-c-1}-\eqref{eq:lip-c-2} and the convexity of $\mathcal{K}_c$ (as in the proof of Proposition~\ref{prop:contraction}), the matrix
\begin{equation*}
    Q = \begin{bmatrix} J_{AA} & J_{A\pi} \\ J_{\pi A} & J_{\pi \pi} \end{bmatrix}
\end{equation*}
satisfies $\tnorm{Q}_1 \le C_c \Delta t$,
while $\tnorm{[J_{\rho A} \ J_{\rho \pi}]}_1 \le C_{\rm adv} + C_c \Delta t$.
To bound the powers of $J_c$, let $\mathbf{x}^{0} = \mathbf{x}$ and $\mathbf{x}^{j+1} = J_c \mathbf{x}^{j}$, and denote with $\mathbf{y}^{j}$ the pair of $A$ and $\pi$ components of $\mathbf{x}^{j}$
and with $x^{j}_\rho$ its $\rho$ component. Since the columns of $J_c$ associated with $\rho$ vanish, $x^{j}_\rho$ does not affect $\mathbf{x}^{j+1}$, and
\begin{equation*}
    \mathbf{y}^{j+1} = Q \, \mathbf{y}^{j}, \qquad x^{j+1}_\rho = [J_{\rho A} \ J_{\rho \pi}] \, \mathbf{y}^{j}.
\end{equation*}
Hence, $\mathbf{y}^{k} = Q^k \mathbf{y}^{0}$ and $x^{k}_\rho = [J_{\rho A} \ J_{\rho \pi}] \, Q^{k-1} \mathbf{y}^{0}$, so that
\begin{equation*}
    \tnorm{J_c^k \mathbf{x}}_1 = \tnorm{\mathbf{y}^{k}}_1 + \tnorm{x^{k}_\rho} \le \left( C_c \Delta t + C_{\rm adv} + C_c \Delta t \right) (C_c \Delta t)^{k-1} \tnorm{\mathbf{x}}_1,
\end{equation*}
which gives \eqref{eq:k-step-c} for $C_c \Delta t \le 1/2$.
Finally, reordering the unknowns as $(A, \pi, \rho)$, $J_c$ becomes block lower triangular with diagonal blocks $Q$ and $0$.
Hence, the nonzero eigenvalues of $J_c$ are those of $Q$, and the spectral radius of $J_c$ is bounded by $\tnorm{Q}_1 \le C_c \Delta t$.
\end{proof}

\begin{remark}\label{rem:comparison}
Propositions~\ref{prop:contraction} and \ref{prop:contraction-coupled} show that the asymptotic convergence rate of both approaches is $\mathcal{O}(\Delta t)$:
by \eqref{eq:k-step-d} and \eqref{eq:k-step-c}, the error after $k$ iterations near the fixed point is bounded by $(C_d \Delta t)^k$ and $(C_{\rm adv} + 1)(C_c \Delta t)^{k-1}$ times the initial error, respectively.
The reason is that both approaches contain the same feedback loops: the loop between $H$ and $\pi$, which is identical in the two approaches, and the loop between $A$ and $\rho$.
In the fully decoupled approach, the latter closes through the lagged reaction term, i.e., $\rho^{k} \to A^{k+1}$ (damped by $\Delta t$) $\to \rho^{k+1}$ (undamped, through the advection term).
In the coupled approach, the Newton linearization removes the lag in the reaction term, but the advection term is lagged instead, so that the loop closes as
$A^{k} \to \rho^{k+1}$ (undamped) $\to A^{k+1}$ (damped by $\Delta t$, through the reaction term in the equation for $A$).
Hence, the two approaches have the same order of convergence in $\Delta t$. However, their constants differ, since the loop between $A$ and $\rho$ closes
through different couplings in the two approaches, and both constants depend on the current state through the bounds in Assumptions~\ref{assump:bounds} and \ref{assump:bounds-coupled}.
As a result, in some states the contraction factor of the coupled approach can be substantially smaller than that of the fully decoupled approach, which can make the difference
between convergence and failure (see Sec.~\ref{sec:num_res}). On the other hand, the coupled approach may pay a transient of about one iteration due to the undamped term $C_{\rm adv}$.
\end{remark}

Because of the term $C_{\rm adv}$, $G_c$ is in general not a contraction in $\tnorm{\cdot}_1$. However, it is a contraction with factor $\mathcal{O}(\Delta t^{1/2})$ in the equivalent norm obtained
by scaling the $A$ component by $\Delta t^{-1/2}$.
Since $\mathcal{K}_c$ is closed, the Banach fixed-point theorem then implies that, for $\Delta t$ sufficiently small, $G_c$ has a unique fixed point in $\mathcal{K}_c$,
i.e., the space- and time-discrete problem at time $t^{n+1}$ has a unique solution in $\mathcal{K}_c$.
Since the least-squares problem may be posed in any norm, the local convergence result in \cite{TothKelley2015} applies.
The result in \cite{Evans2020}, instead, would require the least-squares problem to be posed in this rescaled norm, which differs from $\|\cdot\|_B$ with $B$ given by \eqref{eq:B}.
Since the spectra of $J_d$ and $J_c$ are both contained in a disk of radius $\mathcal{O}(\Delta t)$, we expect AA to provide comparable asymptotic gains for the two approaches (see discussion at the end of Sec.~\ref{sec:AA}),
with possibly one or two additional iterations for the coupled approach due to %
$J_{\rho A}$.

\section{Numerical results}\label{sec:num_res}

The approaches described in Sec.~\ref{sec:decoupled} and
\ref{sec:coupled} were implemented in a solver
based on FEniCSx \cite{BarattaEtal2023,BasixJoss,AlnaesEtal2014} that is
available on GitHub \cite{OurCode}.

We consider the 10 test cases in \cite{Hao2026b}, which show
a rich range of dynamics, including persistent oscillations, splitting and merging of hotspots, and multi-frequency regimes.
All the tests involve
computational domain $\Omega = [0,10]\times [0,10]$ meshed with a structured
grid of $100 \times 100$ elements of side $1/10$. We use first-order Lagrange finite elements on quadrilateral cells
$\mathbb{Q}_1$.
All integrals are evaluated with tensor-product Gauss--Legendre rules whose degree is chosen automatically by the FEniCS Form Compiler (FFCx) from the polynomial degrees of the integrand factors.

We choose stopping
tolerances $\texttt{tol}_1 = \texttt{tol}_2 = \texttt{tol}_3 = 10^{-9}$ {for \eqref{eq:crit}}.
All linear systems, i.e., \eqref{eq:mat1bis}, \eqref{eq:mat2bis}, and \eqref{eq:mat4}, as well as \eqref{eq:mat-coupled} in Sec.~\ref{sec:coupled},
are solved with GMRES \cite{saad2003iterative} with left preconditioning and relative tolerance $10^{-12}$.
As a preconditioner, we use an incomplete LU factorization without fill-in, ILU(0).
At each iteration, GMRES is warm-started, i.e., the solution of the corresponding linear system at the previous iteration is used as {an} initial guess.
In all the tests, we set $\Ast = 1/50$, $A_0 = \Ast + q$, and $\rho_0 = 0.6$. The remaining parameters for each case are summarized in Tab.~\ref{tab:case_summary}.

\begin{table}[htb!]
\caption{Parameter values for the tests from \cite{Hao2026b} considered here.
The regime row indicates whether the case features a stable (S) or oscillatory (O) solution. %
}
\label{tab:case_summary}
\centering
{\small
\setlength{\tabcolsep}{4pt}%
\setlength{\aboverulesep}{0pt}\setlength{\belowrulesep}{0pt}%
\renewcommand{\arraystretch}{1.25}%
\begin{tabular}{l*{10}{!{\color{gray!60}\vrule width 0.4pt}c}}
\toprule
Case & 1 & 2 & 3 & 4 & 5 & 6 & 7 & 8 & 9 & 10 \\
\midrule
Regime  & S     & O     & O      & O       & O      & O      & S      & O      & O      & O \\
\hline
$\eta$  & $0.7$ & $0.3$ & $0.15$ & $0.075$ & $0.15$ & $0.15$ & $0.15$ & $0.15$ & $0.15$ & $0.15$ \\
$q$     & $1.5$ & $1.5$ & $1.5$  & $1.5$   & $0.5$  & $2.5$  & $1.5$  & $1.5$  & $1.5$  & $1.5$ \\
$\tau$  & $5$   & $5$   & $5$    & $5$     & $5$    & $5$    & $0.5$  & $50$   & $5$    & $5$ \\
$\pi_0$ & $0.5$ & $0.5$ & $0.5$  & $0.5$   & $0.5$  & $0.5$  & $0.5$  & $0.5$  & $0.1$  & $1.0$ \\
\bottomrule
\end{tabular}}
\end{table}

In Sec.~\ref{sec:num_nl}, we compare the fully decoupled and coupled approaches without acceleration for all the cases in Tab.~\ref{tab:case_summary},
both for fixed $\Delta t$ and as $\Delta t$ varies.
In Sec.~\ref{sec:num_AA}, we repeat the comparison of the two approaches with AA, and we study how the gain provided by AA depends on the time step and on the depth $m$.
Finally, in Sec.~\ref{sec:num_bounds}, we check numerically if the bounds in Assumptions~\ref{assump:bounds} and \ref{assump:bounds-coupled} hold.

\subsection{Plain fixed-point iterations}
\label{sec:num_nl}

\begin{figure}
\centering
\includegraphics[width=\textwidth]{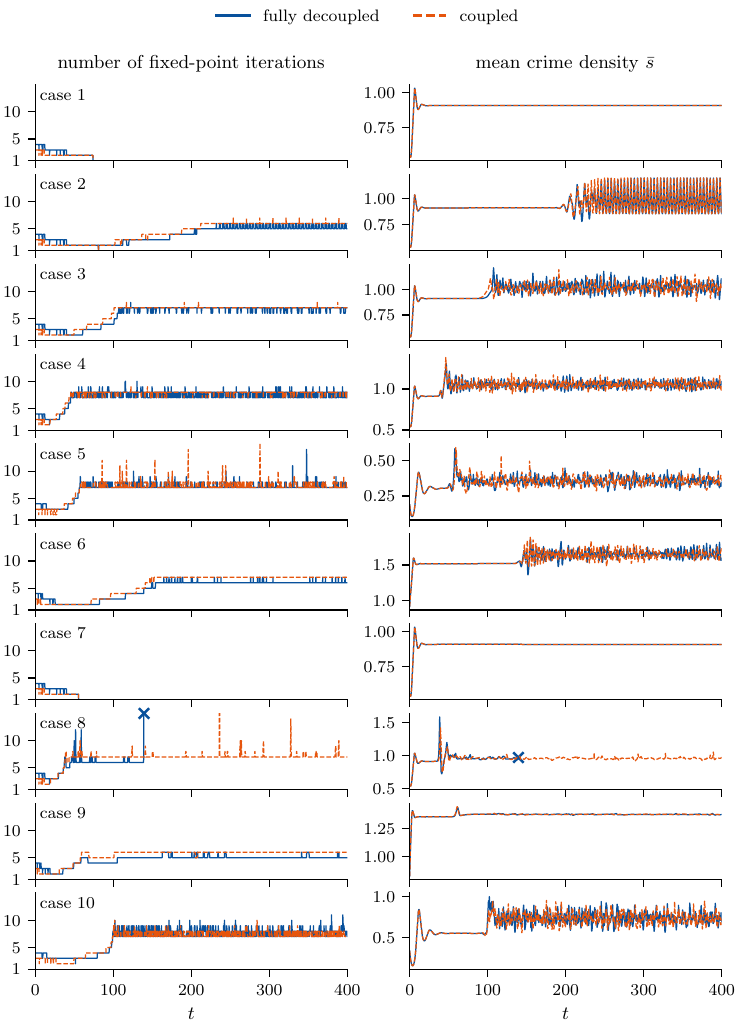}
\caption{Fixed-point iterations without AA - Left: number of fixed-point iterations at each time step. The vertical axis is restricted to $[1, 15]$ and larger values (e.g., in case 8) leave the plotted range. Right: mean crime density $\bar s(t)$.
The simulation of case 8 with the decoupled algorithm terminates at $t = 139.44$, marked by an
$\times$ in both subplots.}
\label{fig:nl_dt002}
\end{figure}

We consider the time interval $[0, 400]$
and set the time step $\Delta t = 0.02$. All simulations reach the end of the time interval, except case 8 using the decoupled algorithm because the algorithm terminates when the fixed-point iteration number reaches 500 without satisfying the stopping criterion ({at} $t = 139.44$).

Let
\[
\bar{s}(t)= \frac{1}{|\Omega|}\int_\Omega \rho A e^{- \pi} \,d\Omega
\]
be the average number of crimes per unit area.
Fig.~\ref{fig:nl_dt002} shows the number of fixed-point iterations without AA
and the computed $\bar{s}(t)$
for the fully decoupled and coupled approaches.
We see that, in the oscillatory regime, the number of fixed-point iterations per time step oscillates in time, mostly following the dynamics of $\bar s(t)$.
To interpret this behavior, note that the number of iterations needed at time step $t^{n+1}$ is approximately%
\begin{equation}\label{eq:estimate_it}
    \frac{\log \left( \texttt{tol} / e^n_0 \right)}{\log \bar\mu^n},
\end{equation}
where \texttt{tol} is the stopping tolerance for the fixed-point algorithm, $e^n_0$ is the initial error, i.e., the distance between the initial iterate and the solution at $t^{n+1}$, and $\bar\mu^n$ is the contraction factor at that time step.
We estimate $\bar\mu^n$ with \eqref{eq:contraction-factor}, taking $k_0 = 1$ and $k_1 = K - 1$, where $K$ is the number of iterations of {the $n$-th time step}. %
The initial error is approximately the change of the solution over one time step. Hence, it tends to oscillate with the dynamics:
it is large when $\bar s(t)$ is rapidly changing and small away from oscillations. In some cases (e.g., cases 2 and 6), we see that the iteration count goes up before $\bar s$ does. This may be due to the formation of hotspots, which affects the solution locally before it is visible in the spatial average $\bar s$.
Tab.~\ref{tab:nl_dt002} reports the mean and maximum number of fixed-point iterations per time step. The coupled approach has the higher mean in 6 of the 10 cases (2, 3, 4, 5, 6, 9). Additionally, each iteration of the coupled scheme requires the solution of a larger system. So, the main advantage of the coupled approach without AA seems to be robustness (no early termination), not efficiency.

\begin{table}[htb!]
\centering
\footnotesize
\setlength{\tabcolsep}{4pt}
\caption{Fixed-point iterations without AA - Mean and maximum number of fixed-point iterations per time step; time mean and standard deviation of $\bar s$ over $t \in [300, 400]$;
time mean $\bar\mu$ and standard deviation over $t \in [300, 400]$ of the average contraction factor $\bar\mu^n$ over the time steps with $K \ge 4$ iterations;
mean over $t \in [300, 400]$ of the contraction in the first iteration $\mu_1$ over the time steps with at least two iterations.
$^*$: computed up to the failure of the fully decoupled approach at $t = 139.44$.}
\label{tab:nl_dt002}
\setlength{\aboverulesep}{0pt}\setlength{\belowrulesep}{0pt}%
\renewcommand{\arraystretch}{1.2}%
\begin{tabular}{l l*{10}{!{\color{gray!60}\vrule width 0.4pt}c}}
\toprule
 & Case & 1 & 2 & 3 & 4 & 5 & 6 & 7 & 8 & 9 & 10 \\
\midrule
\multirow{2}{*}{Iterations, mean} & decoupled & 1.31 & 4.07 & 5.88 & 7.18 & 6.68 & 4.88 & 1.26 & 5.48$^*$ & 4.67 & 6.60 \\
 & coupled & 1.21 & 4.42 & 6.03 & 7.45 & 6.76 & 5.55 & 1.17 & 6.71 & 5.43 & 6.39 \\
\midrule
\multirow{2}{*}{Iterations, max} & decoupled & 4 & 6 & 8 & 10 & 14 & 7 & 4 & 158$^*$ & 6 & 11 \\
 & coupled & 3 & 7 & 8 & 9 & 17 & 7 & 3 & 43 & 6 & 10 \\
\midrule
\multirow{2}{*}{$\bar s$, mean} & decoupled & 0.910 & 0.983 & 1.015 & 1.049 & 0.353 & 1.650 & 0.910 & -- & 1.379 & 0.733 \\
 & coupled & 0.910 & 0.983 & 1.018 & 1.052 & 0.353 & 1.645 & 0.910 & 0.957 & 1.379 & 0.729 \\
\midrule
\multirow{2}{*}{$\bar s$, std.\ dev.} & decoupled & 0.000 & 0.115 & 0.049 & 0.036 & 0.027 & 0.054 & 0.000 & -- & 0.002 & 0.061 \\
 & coupled & 0.000 & 0.115 & 0.049 & 0.038 & 0.025 & 0.045 & 0.000 & 0.014 & 0.002 & 0.062 \\
\midrule
\multirow{2}{*}{$\bar\mu$, mean} & decoupled & -- & 0.021 & 0.044 & 0.071 & 0.059 & 0.034 & -- & -- & 0.024 & 0.074 \\
 & coupled & -- & 0.026 & 0.038 & 0.047 & 0.043 & 0.034 & -- & 0.048 & 0.034 & 0.045 \\
\midrule
\multirow{2}{*}{$\bar\mu$, std.\ dev.} & decoupled & -- & 0.004 & 0.004 & 0.011 & 0.017 & 0.002 & -- & -- & 0.004 & 0.017 \\
 & coupled & -- & 0.002 & 0.002 & 0.004 & 0.007 & 0.002 & -- & 0.013 & 0.003 & 0.007 \\
\midrule
\multirow{2}{*}{$\mu_1$, mean} & decoupled & -- & 0.020 & 0.036 & 0.056 & 0.042 & 0.038 & -- & -- & 0.019 & 0.045 \\
 & coupled & -- & 0.115 & 0.239 & 0.423 & 0.325 & 0.183 & -- & 0.183 & 0.110 & 0.326 \\
\bottomrule
\end{tabular}
\end{table}

From the right column of Fig.~\ref{fig:nl_dt002}, we see that the two approaches give visibly different $\bar s$ curves in cases 3, 4, 5, 6, and 10.
This is expected since the problem is highly nonlinear and the simulation covers a large period of time. However, we point out that the time
means of $\bar s$ in Tab.~\ref{tab:nl_dt002} agree to two or three digits and the standard deviations, also in Tab.~\ref{tab:nl_dt002}, are of the same size.

\begin{figure}[htb!]
\centering
\includegraphics[width=\textwidth]{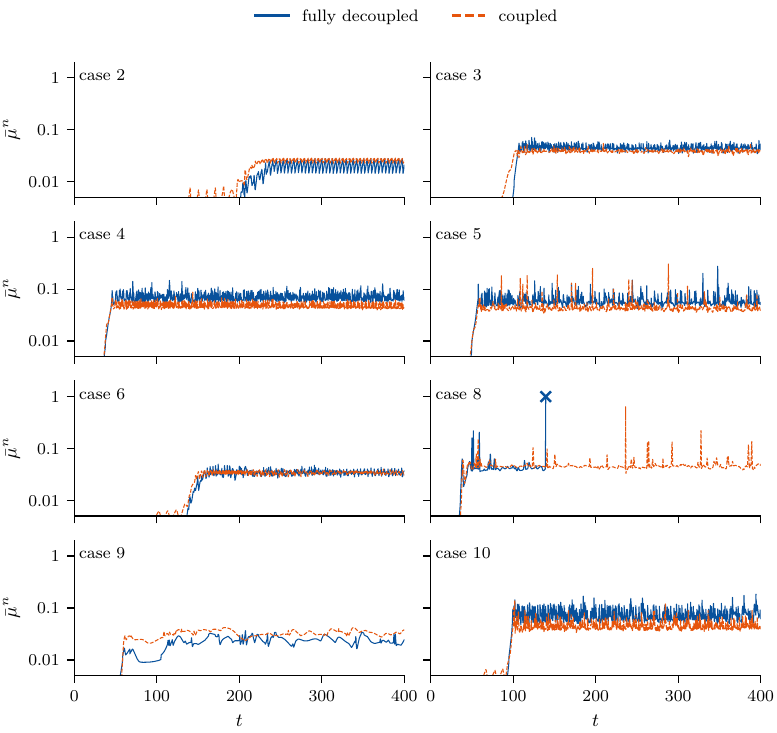}
\caption{Fixed-point iterations without AA - Average contraction factor $\bar\mu^n$ \eqref{eq:contraction-factor} per time step in the Euclidean norm with $k_0 = 1$ and $k_1 = K - 1$, $K$ being the number of iterations of the step. %
Only the cases in the oscillatory regime (see Tab.~\ref{tab:case_summary}) are reported. The time steps with $K \le 3$, those for which $\bar\mu^n$ is undefined (because $\|\mathbf{u}^{K}-\mathbf{u}^{K-1}\|$ is exactly zero) or based on a single ratio, are omitted. The fully decoupled simulation of case 8 terminates at $t = 139.44$, marked by an
$\times$ in the corresponding subplot.}
\label{fig:contraction_dt002}
\end{figure}

To assess the contraction of the fixed-point maps, which is the quantity addressed by the analysis in Sec.~\ref{sec:contraction}
and \ref{sec:contraction-coupled}, Fig.~\ref{fig:contraction_dt002} reports the average contraction factor $\bar\mu^n$ at each time step.
We skip the first iteration ($k_0 = 1$) so that $\bar\mu^n$ is not affected by the transient of the coupled approach
(see Remark~\ref{rem:comparison}), and we compute $\bar\mu^n$ only for time steps with $K \ge 4$.
Although the contraction factor depends on the current state, Fig.~\ref{fig:contraction_dt002} shows that it varies much less than the initial error: in all the oscillatory cases except {for} case 8,
it remains roughly between $0.01$ and $0.1$ throughout the simulation, for both approaches.
Hence, from \eqref{eq:estimate_it}, we deduce that the oscillations in the number of iterations in Fig.~\ref{fig:nl_dt002} are mainly due to the initial error.
Tab.~\ref{tab:nl_dt002}, which reports the time mean of $\bar\mu^n$ and $\mu_1$ for each simulation, confirms the transient predicted by the analysis for the coupled approach:
the contraction factor in the first iteration, $\mu_1$, is between $0.11$ and $0.42$ for the coupled approach and between $0.02$ and $0.06$ for the fully decoupled approach, while the average contraction factor $\bar\mu$, which skips the first iteration, is comparable for the two approaches.
This is consistent with Remark~\ref{rem:comparison}: the undamped term $C_{\rm adv}$ in Proposition~\ref{prop:contraction-coupled} affects the first iteration of the coupled approach,
but not its asymptotic rate.

In case 8, the decoupled algorithm at
$t = 139.44$ settles into a cycle of period two, with $\mathbf{u}^{k+1} = \mathbf{u}^{k-1}$ up to round-off, with relative changes of the order of $10^{-2}$ to $10^{-1}$ that no longer {decrease}. The same happens with a direct linear solver.
To explain this, we computed the dominant eigenvalues of the fixed-point map linearized at the solution of this time step, which we obtained with the coupled approach.
For the fully decoupled approach, the dominant eigenvalue is real and equal to $-1.45$, while all other eigenvalues have modulus at most $0.58$.
Hence, the iteration is not contractive, the error component along the corresponding eigenvector changes sign and grows from one iteration to the next, and the nonlinearity of the map limits it to a stable two-cycle.
For the coupled approach, the dominant eigenvalue at the same solution is $-0.53$, consistent with its observed reduction of the change by a factor of about $0.5$ per iteration. The coupled iteration started from the same state converges in 29 iterations.

\begin{figure}[htb!]
\centering
\includegraphics[width=\textwidth]{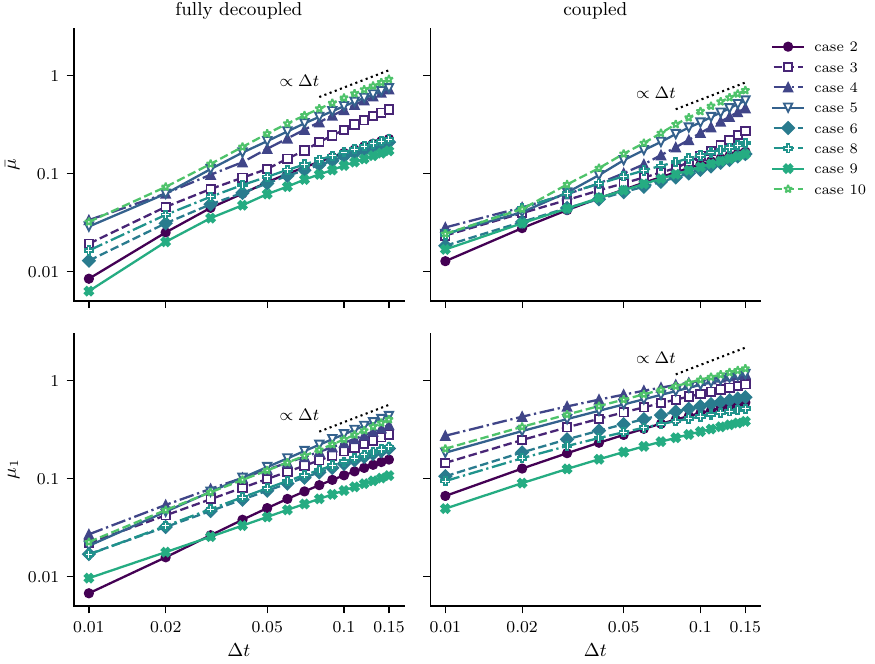}
\caption{Fixed-point iterations without AA for the oscillatory cases -
Top: average contraction factor $\bar\mu$ over iterations $2, \dots, K$ of a single time step as a function of $\Delta t$. Bottom: contraction in the first iteration, $\mu_1 = \|\mathbf{u}^{2} - \mathbf{u}^{1}\| / \|\mathbf{u}^{1} - \mathbf{u}^{0}\|$ as a function of $\Delta t$.
Each step starts from the solution of the coupled approach at $t = 300$, with the solver settings of Tab.~\ref{tab:nl_dt002}.
Left: fully decoupled approach. Right: coupled approach. %
}
\label{fig:contraction_vs_dt}
\end{figure}

To verify the dependence of the contraction factor on $\Delta t$ predicted by Propositions~\ref{prop:contraction} and \ref{prop:contraction-coupled},
we consider a single time step for the oscillatory cases: starting from the solution of the coupled approach at $t = 300$,
we perform one time step with $\Delta t$ ranging from $0.01$ to $0.15$ and compute $\bar\mu$ and $\mu_1$. The results are shown in Fig.~\ref{fig:contraction_vs_dt}.
For the fully decoupled approach, both $\bar\mu$ and $\mu_1$ grow linearly with $\Delta t$, as predicted by Proposition~\ref{prop:contraction}.
Since the bound \eqref{eq:k-step-d} holds from the first iteration, $\mu_1$ and $\bar\mu$ are comparable.
The ratio $\bar\mu / \Delta t$, which estimates the constants $C_d$ and $C_c$, varies roughly between 0.63 and 6.1 for the decoupled scheme and between 1.05 and 4.7 for the coupled scheme, reflecting the dependence of $C_d$ and $C_c$ on the state through the bounds of Assumptions~\ref{assump:bounds} and \ref{assump:bounds-coupled}.
For the coupled approach, $\bar\mu$ also grows approximately like $\Delta t$ and is comparable to that of the fully decoupled approach, as expected from Remark~\ref{rem:comparison}.
However, $\mu_1$ grows more slowly than $\Delta t$. %
This is the transient due to the undamped term $C_{\rm adv}$ in \eqref{eq:k-step-c}, which does not vanish as $\Delta t \to 0$.
Finally, for the largest values of $\Delta t$, $\bar\mu$ reaches about $0.9$ for the fully decoupled approach and $0.7$ for the coupled approach in some cases,
i.e., both plain iterations are close to losing contraction for $\Delta t$ of the order of $0.1$.
The value $\Delta t = 0.02$ used above lies well within the contractive regime, consistently with the values of $\bar\mu^n$ in Fig.~\ref{fig:contraction_dt002}.

\subsection{Anderson-accelerated fixed-point iterations}\label{sec:num_AA}

In this section, we repeat the study of Sec.~\ref{sec:num_nl} with AA($3$) for all the cases in the oscillatory regime in Tab.~\ref{tab:case_summary}. The least-squares problem is posed in norm $\|\cdot\|_B$ with $B$ given by \eqref{eq:B}.
To summarize the results, Tab.~\ref{tab:nl_aa3_dt002} reports the same quantities as in Tab.~\ref{tab:nl_dt002}. In addition, it reports the gain introduced by using AA (measured as the mean $\bar\mu$ without AA in Tab.~\ref{tab:nl_dt002} divided by the mean $\bar\mu$ with AA).

By comparing Tab.~\ref{tab:nl_dt002} and \ref{tab:nl_aa3_dt002}, we see that AA has a limited effect on the average cost of the fixed-point iterations,
but a significant effect on their robustness.
For the fully decoupled approach, the mean number of iterations per time step decreases by at most 10\% (case 3), and it is essentially unchanged in cases 2, 6, and 9.
For the coupled approach, it decreases in cases 3, 4, 6, 8, and 9, and it slightly increases in cases 2, 5, and 10.
This is consistent with the fact that, already without acceleration, both iterations are strongly contractive (see Fig.~\ref{fig:contraction_dt002}).
Moreover, when there are only 4 to 7 iterations per time step, AA($3$) has little room to act: the first iterate is not accelerated, as confirmed by the unchanged values of $\mu_1$,
and the history available to AA remains short.

The benefit of AA is evident in the most demanding cases.
With AA, the maximum number of iterations per time step of the coupled approach decreases from 43 to 14 in case 8 and from 17 to 11 in case 5.
More importantly, with AA the fully decoupled approach completes the simulation of case 8 with at most 15 iterations per time step,
whereas without AA it stalls in a two-cycle at $t = 139.44$ (see Sec.~\ref{sec:num_nl}). We recall that
at that time step, the fixed-point map linearized at the solution has a dominant eigenvalue equal to $-1.45$, while all the other eigenvalues have modulus at most $0.58$.
Hence, the plain iteration is not contractive, but its lack of contraction is due to a single eigenvalue.
This is the situation in which AA is expected to be effective: a few stored iterates suffice to remove the corresponding error component.

For the coupled approach, AA does not improve the asymptotic contraction in cases 4, 5, and 10 (gain $\le 1$), and in cases 5 and 10 the mean number of iterations also increases.
In these cases, we will show below that a change in $\Delta t$ or $m$ might help. We note that in the same cases
the contraction in the first iteration is also the weakest. %
Both with and without AA, $\mu_1$ is 5 to 8 times larger for the coupled approach than for the fully decoupled one,
while the asymptotic contraction factors are comparable. This might explain, at least in part, why the coupled approach needs more iterations on average
in Tab.~\ref{tab:nl_dt002}.

Finally, AA does not alter the computed dynamics in a statistical sense: the time means of $\bar s$ over $[300, 400]$ agree within 1\% with and without AA, for both approaches.
The standard deviations differ more, but this is probably due to the high
nonlinearity of the problem and the long time window.

\begin{table}[htb!]
\centering
\footnotesize
\setlength{\tabcolsep}{4pt}
\caption{Fixed-point iterations with AA(3) for the oscillatory cases - Mean and maximum number of fixed-point iterations per time step; time mean and standard deviation of $\bar s$ over $t \in [300, 400]$;
time mean $\bar\mu$ and standard deviation over $t \in [300, 400]$ of
$\bar\mu^n$ over the time steps with $K \ge 4$ iterations;
mean over $t \in [300, 400]$ of $\mu_1$ %
over the time steps with at least two iterations; gain, i.e., the mean $\bar\mu$ without AA (Tab.~\ref{tab:nl_dt002}) divided by the mean $\bar\mu$ with AA.
$^\dagger$: the fully decoupled approach without AA stops at $t = 139.44$ so the gain cannot be computed. }
\label{tab:nl_aa3_dt002}
\setlength{\aboverulesep}{0pt}\setlength{\belowrulesep}{0pt}%
\renewcommand{\arraystretch}{1.2}%
\begin{tabular}{l l*{8}{!{\color{gray!60}\vrule width 0.4pt}c}}
\toprule
 & Case & 2 & 3 & 4 & 5 & 6 & 8 & 9 & 10 \\
\midrule
\multirow{2}{*}{Iterations, mean} & decoupled & 4.01 & 5.29 & 6.68 & 6.29 & 4.88 & 5.79 & 4.65 & 6.20 \\
 & coupled & 4.45 & 5.79 & 7.35 & 6.83 & 4.97 & 5.83 & 5.23 & 6.65 \\
\midrule
\multirow{2}{*}{Iterations, max} & decoupled & 6 & 7 & 9 & 15 & 6 & 15 & 6 & 9 \\
 & coupled & 6 & 8 & 9 & 11 & 7 & 14 & 6 & 9 \\
\midrule
\multirow{2}{*}{$\bar s$, mean} & decoupled & 0.983 & 1.017 & 1.056 & 0.353 & 1.653 & 0.959 & 1.379 & 0.735 \\
 & coupled & 0.981 & 1.022 & 1.052 & 0.353 & 1.650 & 0.956 & 1.379 & 0.725 \\
\midrule
\multirow{2}{*}{$\bar s$, std.\ dev.} & decoupled & 0.115 & 0.047 & 0.046 & 0.022 & 0.074 & 0.014 & 0.003 & 0.058 \\
 & coupled & 0.115 & 0.056 & 0.035 & 0.029 & 0.088 & 0.015 & 0.002 & 0.059 \\
\midrule
\multirow{2}{*}{$\bar\mu$, mean} & decoupled & 0.012 & 0.028 & 0.048 & 0.042 & 0.020 & 0.019 & 0.014 & 0.052 \\
 & coupled & 0.015 & 0.029 & 0.047 & 0.045 & 0.019 & 0.023 & 0.016 & 0.050 \\
\midrule
\multirow{2}{*}{$\bar\mu$, std.\ dev.} & decoupled & 0.002 & 0.003 & 0.008 & 0.008 & 0.002 & 0.006 & 0.002 & 0.007 \\
 & coupled & 0.002 & 0.004 & 0.007 & 0.008 & 0.002 & 0.006 & 0.001 & 0.007 \\
\midrule
\multirow{2}{*}{$\mu_1$, mean} & decoupled & 0.020 & 0.037 & 0.058 & 0.041 & 0.037 & 0.032 & 0.021 & 0.047 \\
 & coupled & 0.115 & 0.243 & 0.424 & 0.328 & 0.185 & 0.183 & 0.115 & 0.329 \\
\midrule
\multirow{2}{*}{Gain} & decoupled & 1.78 & 1.58 & 1.46 & 1.40 & 1.70 & --$^\dagger$ & 1.72 & 1.42 \\
 & coupled & 1.72 & 1.35 & 1.00 & 0.97 & 1.77 & 2.09 & 2.13 & 0.90 \\
\bottomrule
\end{tabular}
\end{table}

\begin{figure}[htb!]
\centering
\includegraphics[width=\textwidth]{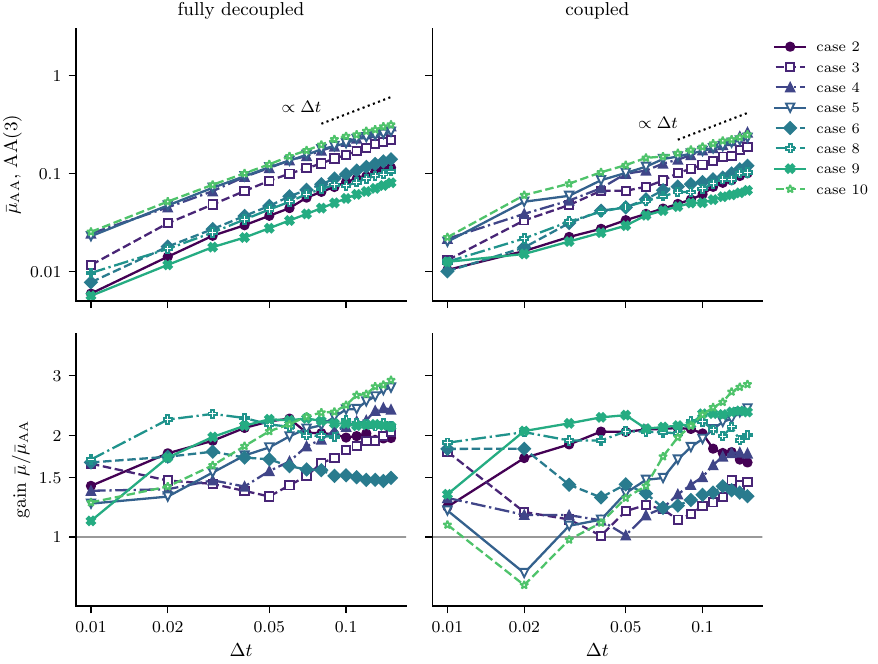}
\caption{Fixed-point iterations with AA(3) for the oscillatory cases -
Top: average contraction factor $\bar\mu_{\rm AA}$ of a single time step, computed as in \eqref{eq:contraction-factor}, as a function of $\Delta t$.
Each step starts from the solution of the coupled approach at $t = 300$, with the solver settings of Tab.~\ref{tab:nl_aa3_dt002}.
Bottom: gain, i.e., $\bar\mu$ without AA (Fig.~\ref{fig:contraction_vs_dt}) divided by $\bar\mu$ with AA, for the same time steps.
Left: fully decoupled approach. Right: coupled approach.}
\label{fig:contraction_vs_dt_aa3}
\end{figure}

Fig.~\ref{fig:contraction_vs_dt_aa3} displays the dependence of the contraction factor on $\Delta t$ when using AA(3). It also shows how the gain varies with the time step. As in the experiment shown in  Fig.~\ref{fig:contraction_vs_dt},
starting from the solution of the coupled approach at $t = 300$,
we perform one time step with $\Delta t$ ranging from $0.01$ to $0.15$ and compute $\bar\mu_{\rm AA}$ and the gain.
Hence, the values for $\Delta t = 0.02$ need not coincide with the values in Tab.~\ref{tab:nl_aa3_dt002}.
For both approaches, the contraction factor with AA still grows approximately like $\Delta t$, i.e., AA does not change the order in $\Delta t$ of the contraction factor but reduces its constant,
in agreement with Remark~\ref{rem:AA-norm}.
All the steps with AA converge, including those with the largest values of $\Delta t$, for which the plain iterations are close to losing contraction (see Fig.~\ref{fig:contraction_vs_dt}).
For the fully decoupled approach, the gain is larger than one for all values of $\Delta t$ and, in most cases, it increases with $\Delta t$.
This is expected: the larger the contraction factor of the plain iteration, the more there is to accelerate.
For the coupled approach, the gain is less regular and varies more among the cases. For $\Delta t = 0.02$, it drops below one in cases 5 and 10, while it is larger than one for $\Delta t = 0.01$ and $\Delta t \ge 0.04$.
Since each value refers to a single time step with few iterations, it is sensitive to the details of the stopping criterion.
Indeed, when averaged over $t \in [300, 400]$, the gain in these cases is close to one (see Tab.~\ref{tab:nl_aa3_dt002}).

\begin{figure}[htb!]
\centering
\includegraphics[width=\textwidth]{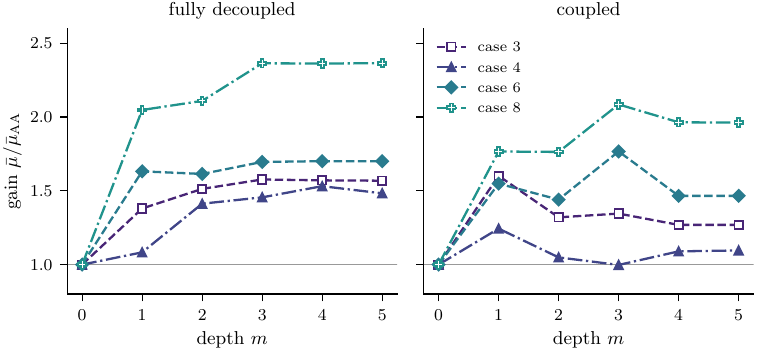}
\caption{Fixed-point iterations with AA($m$) for cases 3, 4, 6 and 8 - Average gain (i.e., $\bar\mu$ without AA divided by $\bar\mu$ with AA($m$)) over the interval [300, 400]
as a function of the depth $m$. %
Left: fully decoupled approach. Right: coupled approach.}
\label{fig:aa_depth_gain}
\end{figure}

Fig.~\ref{fig:aa_depth_gain} shows how the gain depends on the depth $m$ of AA for cases 3, 4, 6, and 8.
We consider the mean over the interval [300, 400].
For case 8 with the fully decoupled approach, which terminates early without AA, the reference run without AA starts from the solution of the coupled approach at $t = 300$.
For the fully decoupled approach, the gain increases with $m$ and saturates at $m = 3$ for cases 3, 6, and 8, while in case 4 it saturates at $m = 4$.
Notice that most of the gain is already obtained with $m = 1$, except in case 4.
For the coupled approach, the gain is generally smaller and depends on $m$ in a less regular way.
In cases 3 and 4, the largest gain is obtained with $m = 1$, and with $m = 3$ the gain is lower. In case 4, AA($3$) brings no gain at all.
In cases 6 and 8, the gain peaks at $m = 3$.
For both approaches, a depth larger than 3 brings at most a small improvement, and only in case 4. This justifies the choice $m = 3$ in Tab.~\ref{tab:nl_aa3_dt002} and Fig.~\ref{fig:contraction_vs_dt_aa3},
although for the coupled approach a smaller depth may be equally effective.

Finally, we comment on the wall time for the tests in this section and the previous.
Tab.~\ref{tab:wall_time} reports the wall time of the fixed-point iterations with and without AA($3$) for case 3 over $t \in [300, 400]$ with $\Delta t = 0.02$, together with a breakdown of the time per iteration. All runs started from the solution of the coupled approach without AA at $t = 300$ and were run serially on one core (Intel Core i9-9900K processor, with the clock frequency fixed at 3.6~GHz).
With our current serial implementation, the assembly of the linear systems dominates the cost: it takes about 85\% of the time per iteration, and its cost is not affected by AA.
The direct cost of AA, i.e., the solution of the least-squares problem and the update, is about 2\% of the time per iteration.
In addition, the linear solves become slightly more expensive with AA, since the average number of GMRES iterations per solve increases.
A possible explanation is that AA removes the last fixed-point iterations of each time step, for which the warm-started linear solves need the fewest GMRES iterations.
Overall, AA makes each iteration about 3\% more expensive, but it reduces the number of iterations by 12\% for the fully decoupled approach and by 5\% for the coupled approach,
so that the wall time decreases by 9\% and 2\%, respectively.

Tab.~\ref{tab:wall_time} also quantifies the different cost per iteration of the two approaches.
An iteration of the coupled approach is about 23\% more expensive than one of the fully decoupled approach, mainly because the solution of the larger linear system
takes about 70\% longer and the assembly about 17\% longer.
Since the two approaches need a similar number of iterations, the coupled approach is %
slower than the fully decoupled one
both without and with AA.
This additional cost is the price for the robustness of the coupled approach without AA, as observed in case 8.
Finally, we note that these timings refer to serial runs. Since the assembly is the dominant cost and it parallelizes well, we expect the relative weight of the linear solves
(and hence of the difference in cost between the two approaches) to increase in parallel runs.

\begin{table}[htb!]
\centering
\footnotesize
\setlength{\tabcolsep}{4pt}
\caption{Wall time of the fixed-point iterations without AA and with AA(3) for case 3 over $t \in [300, 400]$ (5000 time steps): median and, in brackets, minimum and maximum over 5 runs; iterations: total number of fixed-point iterations (identical in all runs);
time per iteration: median wall time divided by the number of iterations, split into assembly, linear solve, AA (least-squares problem and update), and the rest.}
\label{tab:wall_time}
\setlength{\aboverulesep}{0pt}\setlength{\belowrulesep}{0pt}%
\renewcommand{\arraystretch}{1.2}%
\begin{tabular}{l l!{\color{gray!60}\vrule width 0.4pt}c!{\color{gray!60}\vrule width 0.4pt}c!{\color{gray!60}\vrule width 0.4pt}c c c c c}
\toprule
 & & Wall time [s] & Iterations & \multicolumn{5}{c}{Time per iteration [ms]} \\
 & & & & total & assembly & linear solve & AA & rest \\
\midrule
\multirow{2}{*}{Plain} & decoupled & 1674 [1673, 1675] & 34\,250 & 48.87 & 42.70 & 5.63 & -- & 0.53 \\
 & coupled & 2109 [2107, 2123] & 35\,024 & 60.23 & 49.94 & 9.51 & -- & 0.71 \\
\midrule
\multirow{2}{*}{AA(3)} & decoupled & 1522 [1522, 1523] & 30\,184 & 50.43 & 42.79 & 5.94 & 1.08 & 0.62 \\
 & coupled & 2063 [2060, 2080] & 33\,171 & 62.21 & 50.02 & 9.93 & 1.48 & 0.77 \\
\bottomrule
\end{tabular}
\end{table}

Among the methods considered, the fully decoupled approach with AA is the most efficient: in all the oscillatory cases, it needs on average no more iterations per time step than the others, it has the shortest wall time, and it is at least as robust as the others. AA is essential for its robustness: without it, the fully decoupled approach fails in case 8, while with it all the simulations reach the final time, with a maximum number of iterations per time step comparable to that of the coupled approach with AA.

\subsection{Numerical verification of Assumptions~\ref{assump:bounds} and \ref{assump:bounds-coupled}}\label{sec:num_bounds}

We consider the same settings as in Sec.~\ref{sec:num_nl}. To check whether the bounds in Assumptions~\ref{assump:bounds} and \ref{assump:bounds-coupled} hold, Tab.~\ref{tab:bounds} reports the
quantities that appear in both assumptions, computed as extreme values over all the iterates of all the time steps.
For $\rho_h$ and $\pi_h$, it also reports the minima obtained with AA(3) in Sec.~\ref{sec:num_AA}. With AA, $A_h$ and $H_h$ remain positive as well. In the implementation, the accelerated iterate is checked for finite values and, in the coupled approach, for $A_h > 0$. These checks are satisfied in all our simulations.
The minima and maxima of the fields are evaluated at the mesh nodes.
Since the magnitude of the gradient of a $\mathbb{Q}_1$ function on a rectangular cell attains its maximum at one of the corners,
$\max |\nabla A_h|$ and $\max |\nabla H_h|$ are computed from the gradients at the corners of each cell.
Finally, $\rho_h e^{-\pi_h}$, which is not a finite element function, is evaluated at the $3 \times 3$ Gauss points of each cell.
Its maximum should be compared with $1 + 1/\Delta t = 51$ according to the
sufficient condition \eqref{eq:suff_cond} for the matrix in \eqref{eq:mat1bis} to be symmetric positive definite.

\begin{table}[htb!]
\centering
\scriptsize
\setlength{\tabcolsep}{2pt}
\setlength{\aboverulesep}{0pt}\setlength{\belowrulesep}{0pt}
\renewcommand{\arraystretch}{1.2}
\caption{Bounds in Assumptions~\ref{assump:bounds} and \ref{assump:bounds-coupled} with the setting from Sec.~\ref{sec:num_nl}: extreme values over all iterates of all time steps. The iterates of the failing step of the fully decoupled approach in case 8 ($t = 139.44$) are excluded.
Rows labeled AA(3): minima for the runs with AA(3) from Sec.~\ref{sec:num_AA}, taken over the outputs $G(\mathbf{u}^{k})$ of all iterations (oscillatory cases only).
dec.: fully decoupled approach; cou.: coupled approach.}
\label{tab:bounds}
\begin{tabular}{l!{\color{gray!60}\vrule width 0.4pt}l*{10}{!{\color{gray!60}\vrule width 0.4pt}c}}
\toprule
 & Case & 1 & 2 & 3 & 4 & 5 & 6 & 7 & 8 & 9 & 10 \\
\midrule
\multirow{2}{*}{$\min A_h$} & dec. & 0.654 & 0.253 & 0.0581 & 0.0316 & 0.0227 & 0.131 & 0.654 & 0.0324 & 0.149 & 0.0286 \\
 & cou. & 0.654 & 0.258 & 0.0513 & 0.0328 & 0.0236 & 0.198 & 0.654 & 0.0312 & 0.144 & 0.0255 \\
\midrule
\multirow{2}{*}{$\min H_h$} & dec. & 0.550 & 0.313 & 0.0988 & 0.0456 & 0.0109 & 0.341 & 0.541 & 0.287 & 0.0338 & 0.141 \\
 & cou. & 0.550 & 0.318 & 0.116 & 0.0384 & 0.00796 & 0.409 & 0.541 & 0.250 & 0.0326 & 0.100 \\
\midrule
\multirow{4}{*}{$\min \rho_h$} & dec. & 0.600 & 0.360 & 0.0799 & 0.0283 & 0.0214 & 0.137 & 0.600 & 0.0502 & 0.150 & 0.0434 \\
 & cou. & 0.600 & 0.363 & 0.0753 & 0.0268 & 0.0223 & 0.188 & 0.600 & 0.0519 & 0.166 & 0.0469 \\
 & dec., AA(3) & -- & 0.367 & 0.0746 & 0.0246 & 0.0190 & 0.148 & -- & 0.0508 & 0.153 & 0.0476 \\
 & cou., AA(3) & -- & 0.367 & 0.0736 & 0.0249 & 0.0235 & 0.155 & -- & 0.0530 & 0.170 & 0.0572 \\
\midrule
\multirow{4}{*}{$\min \pi_h$} & dec. & 0.500 & 0.0313 & 0.00147 & -0.0325 & -0.00811 & 0.0120 & 0.500 & 0.0421 & $1.3{\cdot}10^{-5}$ & 0.0106 \\
 & cou. & 0.500 & 0.0341 & 0.00339 & -0.0203 & -0.0383 & 0.0150 & 0.500 & 0.0300 & $1.6{\cdot}10^{-5}$ & 0.00246 \\
 & dec., AA(3) & -- & 0.0364 & 0.00252 & -0.0234 & -0.0292 & 0.0136 & -- & 0.0247 & $1.9{\cdot}10^{-5}$ & 0.00461 \\
 & cou., AA(3) & -- & 0.0361 & 0.00306 & -0.0212 & -0.0308 & 0.0135 & -- & 0.0308 & $1.8{\cdot}10^{-5}$ & 0.00469 \\
\midrule[0.8pt]
\multirow{2}{*}{$\max A_h$} & dec. & 1.52 & 4.47 & 9.16 & 12.0 & 12.7 & 9.54 & 1.52 & 12.4 & 7.44 & 8.67 \\
 & cou. & 1.52 & 4.45 & 9.32 & 12.2 & 22.5 & 9.28 & 1.52 & 30.3 & 7.48 & 9.36 \\
\midrule
\multirow{2}{*}{$\max H_h$} & dec. & 0.910 & 2.51 & 3.56 & 4.40 & 6.43 & 4.44 & 1.02 & 1.88 & 8.45 & 2.55 \\
 & cou. & 0.910 & 2.50 & 3.69 & 4.07 & 14.8 & 4.36 & 1.02 & 2.63 & 7.97 & 2.68 \\
\midrule
\multirow{2}{*}{$\max \rho_h$} & dec. & 1.87 & 5.82 & 12.5 & 22.8 & 41.2 & 8.91 & 1.87 & 106 & 4.96 & 22.6 \\
 & cou. & 1.87 & 5.49 & 11.9 & 17.7 & 79.7 & 8.55 & 1.87 & 168 & 4.12 & 27.9 \\
\midrule
\multirow{2}{*}{$\max \pi_h$} & dec. & 0.500 & 1.98 & 3.24 & 3.63 & 5.96 & 3.05 & 0.500 & 1.97 & 1.49 & 5.59 \\
 & cou. & 0.500 & 1.95 & 3.07 & 3.58 & 14.0 & 2.53 & 0.500 & 3.05 & 1.25 & 5.65 \\
\midrule
\multirow{2}{*}{$\max |\nabla A_h|$} & dec. & $2{\cdot}10^{-10}$ & 3.60 & 13.4 & 33.9 & 29.8 & 11.9 & $9{\cdot}10^{-10}$ & 38.5 & 8.36 & 15.9 \\
 & cou. & $2{\cdot}10^{-10}$ & 3.56 & 13.2 & 28.5 & 96.2 & 12.2 & $3{\cdot}10^{-10}$ & 179 & 8.40 & 15.8 \\
\midrule
\multirow{2}{*}{$\max |\nabla H_h|$} & dec. & $9{\cdot}10^{-11}$ & 2.41 & 7.51 & 17.5 & 26.7 & 7.83 & $3{\cdot}10^{-10}$ & 4.52 & 14.5 & 6.54 \\
 & cou. & $1{\cdot}10^{-10}$ & 2.36 & 8.22 & 15.3 & 95.4 & 7.29 & $4{\cdot}10^{-10}$ & 14.9 & 13.7 & 6.51 \\
\midrule[0.8pt]
\multirow{2}{*}{$\max \rho_h e^{-\pi_h}$} & dec. & 1.13 & 3.00 & 4.91 & 9.08 & 14.7 & 3.68 & 1.13 & 74.3 & 2.84 & 6.92 \\
 & cou. & 1.13 & 2.92 & 5.03 & 7.48 & 32.2 & 3.73 & 1.13 & 67.2 & 2.68 & 7.61 \\
\bottomrule
\end{tabular}
\end{table}

From Tab.~\ref{tab:bounds}, we see that,
in cases 4 and 5, $\pi_h$ takes slightly negative nodal values, %
while $A_h$, $H_h$, and $\rho_h$ remain positive in all cases.
The negative values are local: in the converged solutions, they occur at no more than 16 of the 10201 nodes, in 12--14\% of the time steps in case 4 and in less than 1\% in case 5.
They appear with both approaches, with and without AA and with similar magnitude, and they are not removed by reducing $\Delta t$. Hence, they are not caused by the fixed-point iteration, by AA, or by the time discretization, but by the spatial discretization.
The Galerkin discretization of the $\pi$ equation with advection field $2 \nabla H_h / H_h$ does not satisfy a discrete maximum principle when the mesh P\'eclet number, which is proportional to $h \, |\nabla H_h| / H_h$, is large.
This is consistent with cases 4 and 5 having  the smallest diffusion coefficient and the smallest source, respectively, %
and thus exhibiting the sharpest hotspots.  Indeed, they have the largest values of $\max |\nabla H_h|$ and, together with case 9, the smallest values of $\min H_h$ in Tab.~\ref{tab:bounds}. %
Hence, the lower bound $\pi_h \ge 0$ in Assumption~\ref{assump:bounds} is violated on the $100 \times 100$ mesh in these cases, although only by small values, and the fixed-point iterations converge nevertheless.
We expect the bound to hold on sufficiently fine meshes.

Tab.~\ref{tab:bounds} also shows that, in cases 5 and 8, the iterates of the coupled approach reach much larger extreme values than those of the fully decoupled approach: $\max |\nabla A_h| = 96.2$ compared with $29.8$ in case 5, and $179$ against $38.5$ in case 8, with similar differences in $\max A_h$ and $\max \rho_h$.
Since these are extreme values over all iterates, this suggests that the coupled iteration overshoots in states with sharp gradients,
which is where the analysis locates its weakness, i.e., the lagged advection term (constant $C_{\rm adv}$ in Proposition~{\ref{prop:contraction-coupled}}).

The sufficient condition $\rho_h e^{-\pi_h} \le 1 + 1/\Delta t = 51$ for the matrix in \eqref{eq:mat1bis} to be symmetric positive definite is satisfied in all cases except {for} case 8.
In that case, $\max \rho_h e^{-\pi_h}$ reaches $74.3$ for the fully decoupled and $67.2$ for the coupled approach, and in both approaches the maximum is attained by a converged solution in a time step in which the fixed-point iteration struggles:
for the fully decoupled approach at $t = 139.42$ (i.e., right before the failure), which needs 158 iterations,
and for the coupled approach at $t = 236.3$, in the time step with the largest contraction factor (see Fig.~\ref{fig:contraction_dt002}).
Hence, the condition is violated exactly where the fixed-point iteration converges slowly or fails. Since the condition is only sufficient, the matrix may still be positive definite there. However, the iterations struggle when $\rho_h e^{-\pi_h}$ becomes large. %

\section{Conclusions}\label{sec:concl}

We studied the iterative solution of the nonlinear problem arising, at each time step,
from the finite element discretization of a
model for urban crime.
The model couples three PDEs and one ODE describing the interactions among criminal density, environmental attractiveness, delayed crime data, and police deployment.
We considered a fully decoupled fixed-point iteration, which solves for one variable at a time, and proposed a coupled one, in which attractiveness and criminal density are computed together.
We showed that, under suitable bounds on the discrete solution and for sufficiently small $\Delta t$, both iterations converge with an asymptotic rate of order $\Delta t$, and we accelerated them with Anderson acceleration (AA).

The numerical results for ten parameter sets confirm the analysis.
Without acceleration, the average contraction factors are of the order of $10^{-2}$ for $\Delta t = 0.02$ and grow linearly with $\Delta t$.
In one case, the fully decoupled iteration fails because the linearized fixed-point map has an isolated eigenvalue of modulus larger than one, while the coupled iteration converges.
AA removes this failure: with AA, all the simulations reach the final time for both approaches.
Since the plain iterations are already strongly contractive, AA reduces the contraction factor of the fully decoupled approach by about 30\% to 45\%,
but it reduces the number of iterations only slightly, and its benefit for the coupled approach is smaller and less regular.
Overall, the fully decoupled approach with AA is the most efficient: in all the oscillatory cases, it needs on average no more iterations per time step than the other variants, it has the shortest wall time in our timing study, and it is at least as robust as the coupled approach.

\section*{Acknowledgments}
Zhong acknowledges partial support from the National Science Foundation grant AoF-$2225507$. The authors used Claude (Anthropic) to help edit the text and check the proofs. All content was reviewed and verified by the authors, who take full responsibility for it.
\appendix

\section{Proof of Lemma~\ref{lem:nodal}}\label{app:nodal}

\begin{proof}
For each node $\bx_i$, we denote by $\rho_i, A_i, \pi_i$ and $\tilde\rho_i, \tilde A_i, \tilde\pi_i$ the nodal values of the two sets of functions,
with $\delta \rho_i = \rho_i - \tilde\rho_i$, $\delta A_i = A_i - \tilde A_i$, $\delta \pi_i = \pi_i - \tilde\pi_i$,
and with $\delta(\rho A e^{-\pi})_i = \rho_i A_i e^{-\pi_i} - \tilde\rho_i \tilde A_i e^{-\tilde\pi_i}$ the nodal values of $\delta I_h (\rho_h A_h e^{-\pi_h})$.
By Assumption~\ref{assump:bounds}, all the nodal values lie in the box $[0, r_{\max}] \times [a_{\min}, a_{\max}] \times [0, p_{\max}]$, where, since $e^{-\pi} \le 1$ for $\pi \ge 0$,
the first derivatives of $\rho A e^{-\pi}$ with respect to $(\rho, A, \pi)$ are bounded in absolute value by $L$ and the second derivatives by $L_2$.
Since the box is convex, the mean value theorem gives
\begin{equation*}
    |\delta(\rho A e^{-\pi})_i| \le L \left( |\delta \rho_i| + |\delta A_i| + |\delta \pi_i| \right).
\end{equation*}
Moreover, let $\bx_i$ and $\bx_j$ be two nodes of the same element. We write $\delta(\rho A e^{-\pi})_i$ and $\delta(\rho A e^{-\pi})_j$ as integrals of the gradient of $\rho A e^{-\pi}$
along the segments from $(\tilde\rho_i, \tilde A_i, \tilde\pi_i)$ to $(\rho_i, A_i, \pi_i)$ and from $(\tilde\rho_j, \tilde A_j, \tilde\pi_j)$ to $(\rho_j, A_j, \pi_j)$, respectively.
By the gradient bounds in Assumption~\ref{assump:bounds}, corresponding points on the two segments differ by at most $g_{\max} h$ in each component,
so the bound on the second derivatives yields
\begin{align*}
    |\delta(\rho A e^{-\pi})_i - \delta(\rho A e^{-\pi})_j| & \le L \left( |\delta \rho_i - \delta \rho_j| + |\delta A_i - \delta A_j| + |\delta \pi_i - \delta \pi_j| \right) \\
    & \quad + 3 L_2 \, g_{\max} \, h \left( |\delta \rho_j| + |\delta A_j| + |\delta \pi_j| \right).
\end{align*}
On a shape-regular and quasi-uniform mesh, for all $v_h \in V_h(\Omega)$ with nodal values $v_i$, we have the norm equivalences
\begin{equation*}
    \|v_h\|^2 \simeq h^2 \sum_i v_i^2, \qquad \|\nabla v_h\|^2 \simeq \sum_{Q_k \in \mathcal{T}_h} \sum_{\bx_i, \bx_j \in Q_k} (v_i - v_j)^2,
\end{equation*}
with constants that depend only on the shape-regularity and quasi-uniformity constants. The second one follows from a scaling argument,
since in two dimensions the $H^1$ seminorm is invariant under scaling and its kernel on each element consists of the constants.
Applying these equivalences to $\delta I_h (\rho_h A_h e^{-\pi_h})$, $\delta \rho_h$, $\delta A_h$, and $\delta \pi_h$, together with the two bounds above and the fact that each node belongs to a bounded number of elements, gives
\begin{align*}
    \|\delta I_h (\rho_h A_h e^{-\pi_h})\| & \le C_{\rm mesh} L \left( \|\delta \rho_h\| + \|\delta A_h\| + \|\delta \pi_h\| \right), \\
    \|\nabla \delta I_h (\rho_h A_h e^{-\pi_h})\| & \le C_{\rm mesh} L \left( \|\nabla \delta \rho_h\| + \|\nabla \delta A_h\| + \|\nabla \delta \pi_h\| \right) \\
    & \quad + C_{\rm mesh} L_2 \, g_{\max} \left( \|\delta \rho_h\| + \|\delta A_h\| + \|\delta \pi_h\| \right).
\end{align*}
Estimate \eqref{eq:nodal} follows from the definition \eqref{eq:tnorm} of $\tnorm{\cdot}$ and $\Delta t \le 1$.
\end{proof}

\bibliographystyle{plain}
\bibliography{bibliography_verified}
\end{document}